# Multi-block Nonconvex Nonsmooth Proximal ADMM: Convergence and Rates Under Kurdyka–Łojasiewicz Property

Maryam Yashtini[1]

## Abstract

We study the convergence and convergence rates of a multi-block proximal alternating direction method of multipliers (PADMM) for solving linearly constrained separable nonconvex nonsmooth optimization problems. This algorithm is an important variant of the alternating direction method of multipliers (ADMM) which includes a proximal term in each subproblem, to cancel out complicated terms in applications where subproblems are not easy to solve or do not admit a simple closed form solution. We consider an over-relaxation step size in the dual update and provide a detailed proof of the convergence for any step size $\beta \in (0, 2)$. We prove the convergence of the sequence generated by the PADMM by showing that it has a finite length and it is Cauchy. Under the powerful Kurdyka–Łojasiewicz (KŁ) property, we establish the convergence rates for the values and the iterates, and we show that various values of KŁ-exponent associated with the objective function can raise PADMM with three different convergence rates. More precisely, we show that if the (KŁ) exponent $\theta = 0$, the sequence generated by PADMM converges in a finite numbers of iterations. If $\theta \in (0, 1/2]$, then the sequential rate of convergence is $cQ^k$ where $c > 0$, $Q \in (0, 1)$, and $k \in \mathbb{N}$ is the iteration number. If $\theta \in (1/2, 1]$, then $\mathcal{O}(1/k^r)$ rate where $r = (1 - \theta)/(2\theta - 1)$ is achieved.

**Keywords** Nonconvex nonsmooth optimization · Proximal ADMM · Kurdyka–Łojasiewicz property · Convergence · Convergence rates

## 1 Introduction

In this paper, we consider the following nonconvex nonsmooth optimization problem

$$\begin{aligned} \min_{x_1, x_2, \ldots, x_p, y} \quad & \sum_{i=1}^p f_i(x_i) + h(y) \\ \text{s.t.} \quad & \sum_{i=1}^p A_i x_i + By + b = 0, \end{aligned} \quad (1)$$

where $p \geq 1$, $x_i \in \mathbb{R}^{n_i}$ are variables with their coefficient matrix $A_i \in \mathbb{R}^{m \times n_i}$, $i = 1, 2, \ldots, p$, $y \in \mathbb{R}^q$ is the last variable with its coefficient matrix $B \in \mathbb{R}^{m \times q}$, and $b \in \mathbb{R}^m$. The model remains as general without $y$ and $By$; but we keep $y$ and $By$ to simplify the notation. The functions $f_i : \mathbb{R}^{n_i} \to \bar{\mathbb{R}}$ for $i = 1, 2, \ldots, p$ are proper and lower semicontinuous, and $h : \mathbb{R}^q \to \bar{\mathbb{R}}$ is proper smooth function. All functions can be nonconvex. We denote $\mathbf{x} = [x_1, \ldots, x_p] \in \mathbb{R}^n$ where $n = \sum_{i=1}^p n_i$, $\mathbf{A} := [A_1 \ldots A_p] \in \mathbb{R}^{m \times n}$, and $\mathbf{Ax} = \sum_{i=1}^p A_i x_i \in \mathbb{R}^m$. With these notations, the problem (1) can be written

$$\min_{\mathbf{x}, y} \sum_{i=1}^p f_i(x_i) + h(y) \quad \text{s.t.} \quad \mathbf{Ax} + By + b = 0. \quad (2)$$

Such optimization problems arise from a broad spectrum of applications including signal/image processing, sparse optimization, and statistical/machine learning [4,13, 17,34,65,67,68,70].

A method that takes advantage of the special structure of (2) is the alternating direction method of multipliers (ADMM). ADMM is a powerful method which reduces the complexity of the original problem by breaking it into several simpler minimization subproblems. The approximate local solutions to these subproblems are coordinated to find a global solution to the original problem. The ADMM algorithm was proposed by Glowinski et al. [28], motivated by Hestenes [42] and Powell [56], and closely related to the Douglas–Rachford [22] and Peachman–Rachford [55] operator splitting methods. It utilizes the *augmented Lagrangian* associated with (2) given by

$$\begin{gathered} \mathcal{L}^\alpha : \mathbb{R}^n \times \mathbb{R}^q \times \mathbb{R}^m \to \mathbb{R} \\ \mathcal{L}^\alpha(\mathbf{x}, y, z) = \sum_{i=1}^p f_i(x_i) + h(y) + \langle z, \mathbf{Ax} + By + b \rangle + \frac{\alpha}{2} \|\mathbf{Ax} + By + b\|^2, \end{gathered} \quad (3)$$

where $\alpha > 0$ and the vector $z \in \mathbb{R}^m$ is the Lagrange multiplier associated with the constraint $\mathbf{Ax} + By + b = 0$. In each iteration, the augmented Lagrangian is minimized over the primal variables $x_i, i = 1, \ldots, p$ and $y$ separately and one after each other, followed by a dual update for $z$. The iterative scheme of ADMM, known as Gauss–Seidel ADMM [20,39], is outlined below:

$$\begin{aligned} x_i^{k+1} &= \arg\min_{x_i} \mathcal{L}^\alpha(x_{<i}^{k+1}, x_i, x_{>i}^k, y^k, z^k); \quad i = 1, \ldots, p \\ y^{k+1} &= \arg\min_y \mathcal{L}^\alpha(\mathbf{x}^{k+1}, y, z^k); \\ z^{k+1} &= z^k + \alpha\beta(\mathbf{Ax}^{k+1} + By^{k+1} + b). \end{aligned} \quad (4)$$

Here $\beta \in (0, 2)$ is the step size or damping parameter, $x_{<i} := [x_1; \ldots; x_{i-1}] \in \mathbb{R}^{n_1 + \cdots + n_{i-1}}$ and $x_{>i} := [x_{i+1}; \ldots; x_p] \in \mathbb{R}^{n_{i+1} + \cdots + n_p}$ (clearly, $x_{<0}$ and $x_{>p}$ are null variables, which may be used for notational ease). Including the damping parameter $\beta$ in dual update has recently attracted many interests as it can affect the speed of convergence. In [5], the convergence of ADMM with $\beta \in (0, (\sqrt{5} + 1)/2)$ was established. Several works proved the pointwise and ergodic iteration complexity for the PADMM with $\beta \in (0, (\sqrt{5} + 1)/2)$ [18,32] when the functions are convex. In nonconvex setting, the convergence and convergence rate of a linearized PADMM was established in [63] with $\beta \in (0, 2)$.

The two-block ADMM ($p = 1$) and its variants have been extensively studied in terms of convergence analysis and numerical implementation to solve convex problems—see [5,10,17,19,21,24,29–31,34,35,40,41,43,64]. The convergence of multi-block ADMM ($p \geq 2$) requires stronger assumptions than two-block ADMM [11,14,15,36,51,53]. For example, in [36] the authors show that all involved functions need to be strongly convex and the $\alpha > 0$ parameter needs to be restricted. In [14], the existence of at least two orthogonal coefficient matrices are essential to ensure the convergence. Lin et al. [53] proved the convergence of the three-block ADMM ($p = 2$) for solving (2) under the assumptions that $h$ is strongly convex, $\nabla h$ has Lipschitz continuous, $f_1$ and $f_2$ are convex, $A_1$ and $A_2$ are full column rank, and $B$ is an identity matrix. Some researchers modified the standard multi-block ADMM to guarantee convergence without imposing strong assumptions. For instance, in [38,39] the authors proposed some contractive substitution into ADMM and proved that the new framework converges without strong convexity assumption. In [39], the authors combined the ADMM with either forward or backward substitution procedure and they showed that convergence of this framework can be easily proved from contraction perspective, and its local linear convergence rate is provable if certain standard error bound condition is assumed.

In [20], the authors showed that the multi-block Jacobi ADMM converges in the following two cases: (i) all coefficient matrices are mutually near orthogonal and have full column rank, or (ii) proximal terms are added to the subproblems, and they proved $o(1/k)$ convergence rate. In [37], the authors combined the Jacobi ADMM with a relaxation step and proved the convergence and $O(1/k)$ rate in both ergodic and nonergodic senses. In [43] the authors shows that when the step size $\beta$ is made sufficiently small the sum of the primal and dual optimality gaps decreases after each iteration.

The ADMM algorithm generally fails to solve nonconvex (possibly nonconvex) optimization problems; however, its great performance on some practical applications [46,49,58,61,65,66] has encouraged researchers to study underlying conditions for which nonconvex nonsmooth ADMM converges. Wang and Yin [60] established the convergence of multi-block ADMM with $\beta = 1$ in the nonconvex nonsmooth setting, for both separable and nonseparable objective functions. For the separable objective function, they assume that $f_2, \ldots, f_p$ satisfy prox-regularity. When the functions are either smooth and nonconvex or convex and nonsmooth, Hong et al. [44] studied the convergence of multi-block ADMM for solving a family of nonconvex consensus and sharing problems in which all coefficient matrices are assumed to be full rank. Wang et al. [59] studied the so-called Bregman ADMM and include the standard ADMM as a special case. They established the convergence under the assumptions where $f_i$'s are strongly convex, there exists a $\beta_0 > 0$ such that $h - \beta_0 \|\nabla h\|^2$ is lower bounded, and $B$ is full row rank.

In [54], the authors studied the convergence of a full linearized ADMM. The work [45] considered a structured nonconvex nonsmooth optimization, different from (2), and showed that the proposed framework converges. The authors analyzed its iteration complexity. In [33], Guo et al. considered a two-block ADMM to minimize a sum of two nonconvex functions with linear constraints and they proved the convergence by assuming that the generated sequence is bounded. Yang et al. [62] considered a three-block ADMM to solve a special class of nonconvex and nonsmooth problem with applications to background/foreground extraction in which $A_1^T A_1$ and $A_2^T A_2$ are assumed to be positive definite and $B = -I$. Li and Pong [47] proved the global convergence of two-block ADMM under assumption that $h$ is twice differentiable and has bounded Hessian.

In this paper, we consider a variant/modification of ADMM that is called Proximal ADMM (PADMM). To state the algorithm let $\{Q_i\}_{i=1,\ldots,p} \subseteq \mathbb{R}^{n_i \times n_i}$ and $P \in \mathbb{R}^{q \times q}$ be symmetric positive semidefinite matrices and $\|x\|_M^2 = x^T M x$, then the nonconvex possibly nonsmooth PADMM algorithm for solving (2) is given by:

$$\begin{aligned}
x_i^{k+1} &= \arg\min_{x_i} \mathcal{L}^\alpha(x_{<i}^{k+1}, x_i, x_{>i}^k, y^k, z^k) + \tfrac{1}{2}\|x_i - x_i^k\|_{Q_i}^2; \quad i = 1, \ldots, p \\
y^{k+1} &= \arg\min_y \mathcal{L}^\alpha(\mathbf{x}^{k+1}, y, z^k) + \tfrac{1}{2}\|y - y^k\|_P^2; \\
z^{k+1} &= z^k + \alpha\beta(\mathbf{A}x^{k+1} + By^{k+1} + b).
\end{aligned} \quad (5)$$

This variant of the ADMM algorithm is important, especially in the applications in which it is expensive to solve some or all subproblems exactly [12,17,18,34,35,39,50, 52]. When the subproblems are not easy to solve a suitable proximal terms such as $\frac{1}{2}\|x_i - x_i^k\|_{Q_i}^2$ and $\frac{1}{2}\|y - y^k\|_P^2$ with proper choices of $Q_i, i = 1, \ldots p$ and $P$ are often added to the standard ADMM to cancel out complicated terms [23,26]. For example, by setting $Q_i = \frac{\alpha}{\tau'} I_{n_i} - \alpha A_i^T A_i$ with $\tau' > 0$ the $x_i$'s subproblems in the PADMM algorithm (5) reduce to prox-linear subproblems [16]. Prox-linear subproblems are easier to compute in various applications in signal and image processing, statistics, machine learning, etc. For instance, if $f_i$ is $\ell_1$ norm, solution is given in the closed form by so-called soft thresholding, if $f_i$ is nuclear norm, then singular value soft thresholding is used. If the function $h$ is quadratic with a positive semidefinite hessian $H(h) = \nabla^2 h(x)$, letting $P = \frac{1}{s} I_q - H(h) - \alpha B^T B$ with $s > 0$ gives rise to a gradient descent step $y^{k+1} = y^k - sg^k$ where $g^k = \nabla h(x^k) + \alpha B^T (\mathbf{A}x^{k+1} + By^k + b + z^k/\alpha)$.

The main contribution of this paper is to establish the convergence and convergence rate of the PADMM Algorithm (5) when the objective function is nonconvex and possibly nonsmooth. Our algorithm and theoretical analysis differ from the existing works [44,54,59,60]; they prove the convergence of the sequence, while we prove the convergence and also the convergence rates. In Theorem 1, we prove that the sequence generated by the PADMM has a finite length, it is Cauchy; hence, it is convergent to a stationary point of (2). The function and sequential convergence rates are established under the assumption that the modified augmented Lagrangian satisfies the KŁ property (see Sect. 2). The KŁ property has been used for analyzing the convergence rate of several first-order methods [2,27]. The KŁ property and its associated KŁ exponent have their roots in algebraic geometry; and they describe a qualitative relationship between the value of a suitable potential function (depending on the optimization model and the algorithm being considered) and some first-order information (gradient or subgradient) of the potential function. We will discuss three possible rates based on the KŁ parameter.

The rest of this paper is organized as follows. Section 2 provides some preliminary results and notations. In Sect. 3, we drive some fundamental properties for the sequence generated by the PADMM algorithm (5) based on the augmented Lagrangian function. In Sects. 4 and 5, we prove the global convergence and convergence rates, respectively. Section 6 provides some practical applications and insight on how to choose the proximal terms. We end the paper with some concluding remarks and future directions in Sect. 7.

## 2 Notation and Preliminaries

Throughout this paper, we denote $\mathbb{R}$ as the real number set while $\mathbb{Z}$ as the set of integers. The set $\bar{\mathbb{R}} := \mathbb{R} \cup \{+\infty\}$ is the extended real numbers, $\mathbb{R}_+$ is the positive real number set, and $\mathbb{Z}_+$ is the set of positive integers. Given the matrix $X$, $\text{Im}(X)$ denotes its image. We denote by $I_n$ the $n \times n$ identity matrix. The minimum and maximum eigenvalues of the matrix $X \in \mathbb{R}^{n \times n}$ are denoted, respectively, by $\lambda_{\min}(X)$, and $\lambda_{\max}(X)$, and $\lambda_{++}^X$ refers to the smallest strictly positive eigenvalue of $X$. The Euclidean scalar product of $\mathbb{R}^n$ and its corresponding norms are, respectively, denoted by $\langle \cdot, \cdot \rangle$ and $\|\cdot\| = \sqrt{\langle \cdot, \cdot \rangle}$. For a vector $x \in \mathbb{R}^n$, we use $\|x\|_1$ to denote $\ell_1$ norm and $\|x\|_0$ to denote the number of entries in $x$ that are nonzero ("$\ell_0$ norm"). For a (nonempty) closed set $D \subset \mathbb{R}^n$, the indicator function $\delta_D$ is defined as

$$\delta_D(x) := \begin{cases} 0 & \text{if } x \in D \\ \infty & \text{otherwise.} \end{cases}$$

If $n_1, \ldots, n_p \in \mathbb{Z}_+$ and $p \in \mathbb{Z}_+$, then for any $v := (v_1, \ldots, v_p) \in \mathbb{R}^{n_1} \times \mathbb{R}^{n_2} \times \cdots \times \mathbb{R}^{n_p}$ and $v' := (v'_1, \ldots, v'_p) \in \mathbb{R}^{n_1} \times \mathbb{R}^{n_2} \times \cdots \times \mathbb{R}^{n_p}$ the Cartesian product and its norm are defined by

$$\ll v, v' \gg = \sum_{i=1}^{p} \langle v_i, v'_i \rangle \quad \frac{1}{\sqrt{p}} \sum_{i=1}^{p} \|v_i\| \leq |||v||| = \sqrt{\sum_{i=1}^{p} \|v_i\|^2} \leq \sum_{i=1}^{p} \|v_i\|.$$

Let $\Phi : \mathbb{R}^d \to \bar{\mathbb{R}}$ be a proper and lower semicontinuous function. The domain of $\Phi$, denoted dom $\Phi$, is defined by dom $\Phi := \{x \in \mathbb{R}^d : \Phi(x) < +\infty\}$. For any $x \in \text{dom } \Phi$, the *Fréchet (viscosity) subdifferential* of $\Phi$ at $x$, denoted $\hat{\partial}\Phi(x)$, is defined by

$$\hat{\partial}\Phi(x) = \left\{ s \in \mathbb{R}^d : \liminf_{\substack{y \neq x \\ y \to x}} \frac{\Phi(y) - \Phi(x) - \langle s, y - x \rangle}{\|y - x\|} \geq 0 \right\}.$$

For $x \notin \text{dom } \Phi$, then $\partial\hat{\Phi}(x) = \emptyset$. The *limiting (Mordukhovich) subdifferential*, or simply the subdifferential for short, of $\Phi$ at $x \in \text{dom } \Phi$, denoted $\partial\Phi(x)$, is defined by

$$\partial\Phi(x) := \{s \in \mathbb{R}^d : \exists x^k \to x, \ \Phi(x^k) \to \Phi(x) \text{ and } s^k \in \hat{\Phi}(x^k) \to s \text{ as } k \to +\infty\}.$$

For any $x \in \mathbb{R}^d$, the above definition implies $\partial\hat{\Phi}(x) \subset \partial\Phi(x)$, where the first set is convex and closed while the second one is closed ([57], Theorem. 8.6).

Let $\{x^k\}_{k \in \mathbb{N}}$ and $\{s^k\}_{k \in \mathbb{N}}$ be sequences such that $x^k \to x^*$, $s^k \to s^*$, $\Phi(x^k) \to \Phi(x^*)$ and $s^k \in \partial\Phi(x^k)$ as $k \to \infty$, then $s^* \in \partial\Phi(x^*)$. The well-known Fermat's rule "$x \in \mathbb{R}^d$ is a local minimizer of $\Phi$, then $\partial\Phi(x) \ni 0$" remains unchanged. If $x \in \mathbb{R}^d$ such that $\partial\Phi(x) \ni 0$ the point $x$ is called a critical point. We denote by crit $\Phi$ the set of *critical points* of $\Phi$, that is crit $\Phi = \{x \in \mathbb{R}^d : 0 \in \partial\Phi(x)\}$. When $\Phi$ is convex the two sets coincide and

$$\partial\hat{\Phi}(x) = \partial\Phi(x) = \{s \in \mathbb{R}^d : \Phi(y) \geq \Phi(x) + \langle s, y - x \rangle \ \forall y \in \mathbb{R}^d\}.$$

Let $\Omega$ be a subset of $\mathbb{R}^d$ and $x$ is any point in $\mathbb{R}^d$. The distance from $x$ to $\Omega$, denoted $\text{dist}(x, \Omega)$, is defined by $\text{dist}(x, \Omega) = \inf\{\|x - z\| : z \in \Omega\}$, and the set of points in $\Omega$ that achieve this infimum (the projection of $x$ onto $\Omega$) is defined by $\text{Proj}_\Omega(x)$. The set $\text{Proj}_\Omega(x)$ becomes a singleton if $\Omega$ is a closed convex set. If $\Omega = \emptyset$, then $\text{dist}(x, \Omega) = +\infty$ for all $x \in \mathbb{R}^d$. For any real-valued function $\Phi$ on $\mathbb{R}^d$, we have

$$\text{dist}(0, \partial\Phi(x)) = \inf\{\|s^*\| : s^* \in \partial\Phi(x)\}.$$

Let $F : \mathbb{R}^n \times \mathbb{R}^m \to \bar{\mathbb{R}}$ be a lower semicontinuous function. The subdifferentiation of $F$ at the point $(\hat{x}, \hat{y})$ is defined by $\partial F(\hat{x}, \hat{y}) = (\partial_x F(\hat{x}, \hat{y}), \partial_y F(\hat{x}, \hat{y}))$, where $\partial_x F$ and $\partial_y F$ are, respectively, the differential of the function $F(\cdot, y)$ when $y \in \mathbb{R}^m$ is fixed, and $F(x, \cdot)$ when $x$ is fixed.

**Lemma 1** *Let $\Phi : \mathbb{R}^d \to \mathbb{R}$ be Fréchet differentiable such that its gradient is Lipschitz continuous with constant $L_\Phi > 0$. Then for every $u, v \in \mathbb{R}^d$ and every $\xi \in [u, v] = \{(1-t)u + tv : t \in [0, 1]\}$ it holds*

$$\Phi(v) \leq \Phi(u) + \langle \nabla\Phi(\xi), v - u \rangle + \frac{L_\Phi}{2}\|v - u\|^2. \tag{6}$$

**Proof** Let $u, v \in \mathbb{R}^d$, and $\xi = (1-t)u + tv$ for $t \in [0, 1]$. Then we have

$$\Phi(v) - \Phi(u) = \int_0^1 \langle \nabla\Phi((1-r)u + rv), v - u \rangle \mathrm{d}r$$
$$= \int_0^1 \langle \nabla\Phi((1-r)u + rv) - \nabla\Phi(\xi), v - u \rangle \mathrm{d}r + \langle \nabla\Phi(\xi), v - u \rangle.$$

Since $\Phi$ has Lipschitz continuous gradients with constant $L_\Phi > 0$, then we have

$$|\Phi(v) - \Phi(u) - \langle \nabla\Phi(\xi), v - u \rangle| = \left|\int_0^1 \langle \nabla\Phi((1-r)u + rv) - \nabla\Phi(\xi), v - u \rangle \mathrm{d}r\right|$$

$$\leq \int_0^1 \left\| \nabla \Phi((1-r)u + rv) - \nabla \Phi(\xi) \right\| \cdot \left\| v - u \right\| \mathrm{d}r \leq L_\Phi \|v - u\|^2 \int_0^1 |r - t| \mathrm{d}r$$

$$= L_\Phi \|v - u\|^2 \left( \int_0^t (-r + t)\mathrm{d}r + \int_t^1 (r - t)\mathrm{d}r \right)$$

$$= L_\Phi \left( \tfrac{1}{2} - t(1-t) \right) \|v - u\|^2 \leq \tfrac{L_\Phi}{2} \|v - u\|^2.$$

□

Let $\Phi : \mathbb{R}^d \to \bar{\mathbb{R}}$ be a proper lower semicontinuous function. For $-\infty < \eta_1 < \eta_2 \leq +\infty$, we define $[\eta_1 < \Phi < \eta_2] = \{x \in \mathbb{R}^d : \eta_1 < \Phi(x) < \eta_2\}$. Let $\eta \in (0, +\infty]$. We denote by $\Psi_\eta$ the set of all continuous concave functions $\psi : [0, \eta] \to [0, +\infty)$ that $\psi(0) = 0$ and $\psi$ is continuously differentiable on $(0, \eta)$ with $\psi'(s) > 0$ over $(0, \eta)$.

**Definition 1** Let $\Phi : \mathbb{R}^d \to \bar{\mathbb{R}}$ be a proper lower semicontinuous function. We say that $\Phi$ has the Kurdyka–Łojasiewicz (KŁ) property at $x^* \in \mathrm{dom}\, \partial \Phi$ if there exists a neighborhood $U$ of $x^*$, $\eta \in (0, \infty)$ and a continuous concave function $\psi \in \Psi_\eta$ such that

$$\psi'\bigl(\Phi(x) - \Phi(x^*)\bigr) \mathrm{dist}\bigl(0, \partial \Phi(x)\bigr) \geq 1. \tag{7}$$

If $\Phi$ satisfies the property at each point of $\mathrm{dom}\, \partial \Phi$, then $\Phi$ is a KŁ function.

**Definition 2** Let $\Phi : \mathbb{R}^d \to \bar{\mathbb{R}}$ be a proper lower semicontinuous function that takes constant value on $\Omega$ and satisfies the KŁ property at each point of $\Omega$. We say $\Phi$ satisfies a KŁ property on $\Omega$ if there exists $\epsilon > 0$, $\eta > 0$, and $\psi \in \Psi_\eta$ such that for every $x^* \in \Omega$ and every element $x$ belongs to the intersection $\{x \in \mathbb{R}^d : \mathrm{dist}(x, \Omega) < \epsilon\} \cap [\Phi(x^*) < \Phi(x) < \Phi(x^*) + \eta]$, (7) holds.

**Definition 3** Let the proper lower semicontinuous function $\Phi$ satisfying the KŁ property at $x^* \in \mathrm{dom}\partial \Phi$, and the corresponding function $\psi \in \Psi_\eta$ can be chosen as $\psi(s) = \bar{c} s^{1-\theta}$ for some $\bar{c} > 0$ and $\theta \in [0, 1)$, i.e., there exist $c > 0$ and a neighborhood $U$ of $x^*$ and $\eta \in (0, +\infty]$ such that

$$(\Phi(x) - \Phi(x^*))^\theta \leq c\, \mathrm{dist}(0, \partial \Phi(x)), \quad \forall x \in U \tag{8}$$

holds. Then we say $\Phi$ has the KŁ property at $x^*$ with an exponent $\theta$. This asserts that $(\Phi(x) - \Phi(x^*))^\theta / \mathrm{dist}(0, \partial \Phi(x))$ remains bounded around $x^*$.

This definition encompasses broad classes of functions arising in practical optimization problems. For example, if $f$ is a proper closed semialgebraic function, then $f$ is a KŁ function with exponent $\theta \in [0, 1)$ [3]. The function $l(Ax)$, where $l$ is strongly convex on any compact set and it is twice differentiable, and $A \in \mathbb{R}^{m \times n}$, is a KŁ function. Convex piecewise linear quadratic function, $\|x\|_1$, $\|x\|_0$, $\gamma \sum_{i=1}^k |x_{[i]}|$ where $|x_{[i]}|$ is the $i$th largest (in magnitude) entry in $x$, $k \leq n$ and $\gamma \in (0, 1]$, $\delta_\Delta(x)$ where $\Delta = \{x \in \mathbb{R}^n : e^T x = 1, x \geq 0\}$, and least square problems with SCAD [25] or MCP [69] regularized functions are all KŁ functions. The KŁ property characterizes the local geometry of a function around the set of critical points.

If $\Phi$ is continuously differentiable then $\partial \Phi(x) = \nabla \Phi(x)$, the inequality (8) becomes $\Phi(x) - \Phi(x^*) \leq c\|\nabla \Phi(x)\|^{\frac{1}{\theta}}$, which generalizes the Polyak–Łojasiewics condition $\Phi(x) - \Phi(x^*) \leq \mathcal{O}(\|\nabla \Phi(x)\|^2)$ ($\theta = \frac{1}{2}$). We refer interested readers to [1–3,6–9,27] for more properties of KŁ functions and illustrating examples.

## 3 Augmented Lagrangian-Based properties

In this section, we establish some properties for the PADMM algorithm (5). We begin by making the following general assumptions on the optimization problem (1):

A1. $f_i, i = 1, \ldots, p$ are lower semicontinuous, and coercive that is

$$\lim_{\|x_i\| \to \infty} f_i(x_i) = +\infty;$$

A2. $h$ is bounded from below, and is $L_h$ Lipschitz differentiable that is

$$\|\nabla h(v) - \nabla h(u)\| \leq L_h \|v - u\|, \quad \forall u, v \in \mathbb{R}^q;$$

A3. $\mathrm{Im}(\mathbf{A}) \subseteq \mathrm{Im}(B)$ and $b \in \mathrm{Im}(B)$;
A4. For $i = 1, \ldots, p$ one of the followings hold:

  (i) $\exists q_i > 0$ such that $Q_i \succeq q_i I_{n_i}$; or
  (ii) $\exists \epsilon_i > 0$ such that $f_i(v) - f_i(\tilde{v}) - \langle \tilde{s}_i, v - \tilde{v} \rangle \geq -\frac{\epsilon_i}{2}\|v - \tilde{v}\|^2$ for $v, \tilde{v} \in \mathbb{R}^{n_i}$ and $\tilde{s}_i \in \partial f_i(\tilde{v})$. This is equivalent to the weak convexity of the function $f_i$, $i = 1, \ldots, p$.

Throughout this paper, we frequently refer to some quantities which for simplicity they are summarized as follow:

$$\begin{aligned}
&\rho(\beta) := 1 - |1 - \beta|, \\
&c_1 := \frac{1}{\alpha \rho(\beta) \lambda_{++}^{B^T B}}, & c_2 := \frac{2\beta}{\alpha \rho(\beta)^2 \lambda_{++}^{B^T B}}, \\
&c_3 := c_2 \|P\|^2, & c_4 := c_2 (\|P\| + L_h)^2, \\
&c_5 := \beta^{-1}|1 - \beta|c_1, & c_6 := \frac{|1-\beta|}{2\alpha \beta^2 \lambda_{++}^{B^T B}}, \\
&D = 2P + \alpha B^T B - L_h I_m, & \bar{D} = D - 2\epsilon_0(c_3 + c_4)I,
\end{aligned} \quad (9)$$

where $\epsilon_0 > 1$, $\lambda_{++}^{B^T B}$ is the smallest strictly positive eigenvalue of $B^T B$, and $I_m$ denotes a $m \times m$ identity matrix.

**Lemma 2** (Subgradient bound) *Let $\{(\mathbf{x}^k, y^k, z^k)\}_{k \geq 0}$ be a sequence generated by the PADMM algorithm (5). Then there exists $d^k := (\{d_{x_i}^k\}_{i=1}^p, d_y^k, d_z^k) \in \partial \mathcal{L}^\alpha(\mathbf{x}^k, y^k, z^k)$ and a constant $\rho > 0$ such that*

$$|||d^k||| \leq \rho \Big( \sum_{i=1}^p \|\Delta x_i^k\| + \|\Delta y^k\| + \|\Delta z^k\| \Big), \quad (10)$$

*where*

$$\rho := \max \left\{ \|Q\| + \|A\|^2, \quad \|P\| + \alpha \|A\| \|B\|, \quad \|A\| + \|B\| + \frac{1}{\alpha\beta} \right\}, \quad (11)$$

$\|Q\| := \sum_{i=1}^{p} \|Q_i\|$ and $\|A\| := \sum_{i=1}^{p} \|A_i\|$, and for any sequence $\{u^k\}_{k \geq 0}$, $\Delta u^{k+1} = u^{k+1} - u^k$.

**Proof** Let $k \geq 1$ be fixed. By the optimality condition of the $x_i$ subproblem of PADMM algorithm for $i = 1, \ldots, p$, we have

$$-A_i^T z^{k-1} - \alpha A_i^T \left( A_{<i} x_{<i}^k + A_i x_i^k + A_{>i} x_{>i}^{k-1} + B y^{k-1} + b \right) - Q_i \Delta x_i^k \in \partial f_i(x_i^k).$$

Taking partial differential of $\mathcal{L}^\alpha$ with respect to $x_i$ and evaluating the result at the point $(\mathbf{x}^k, y^k, z^k)$ yield

$$\partial_{x_i} \mathcal{L}^\alpha(\mathbf{x}^k, y^k, z^k) = \partial f_i(x_i^k) + A_i^T z^k + \alpha A_i^T \left( A\mathbf{x}^k + By^k + b \right).$$

Therefore,

$$d_{x_i}^k := A_i^T \Delta z^k + \alpha A_i^T A_{>i} \Delta x_{>i}^k + \alpha A_i^T B \Delta y^k - Q_i \Delta x_i^k \in \partial_{x_i} \mathcal{L}^\alpha(\mathbf{x}^k, y^k, z^k). \quad (12)$$

The optimality criterion of $y$ subproblem gives

$$\nabla h(y^k) = -B^T z^{k-1} - \alpha B^T (A\mathbf{x}^k + By^k + b) - P \Delta y^k.$$

By taking partial differential of $\mathcal{L}^\alpha$ with respect to $y$ and evaluating it at $(\mathbf{x}^k, y^k, z^k)$, we obtain

$$\nabla_y \mathcal{L}^\alpha(\mathbf{x}^k, y^k, z^k) = \nabla h(y^k) + B^T z^k + \alpha B^T \left( A\mathbf{x}^k + By^k + b \right).$$

The last two equations follow that

$$d_y^k := B^T \Delta z^k - P \Delta y^k = \nabla_y \mathcal{L}^\alpha(\mathbf{x}^k, y^k, z^k). \quad (13)$$

Lastly, by the $z$ subproblem we have $A\mathbf{x}^k + By^k + b = \frac{1}{\alpha\beta} \Delta z^k$; hence,

$$d_z^k := \frac{1}{\alpha\beta} \Delta z^k = \nabla_z \mathcal{L}^\alpha(\mathbf{x}^k, y^k, z^k). \quad (14)$$

Hence, by (12), (13), and (14), we have $d^k := \left( \{d_{x_i}^k\}_{i=1}^p, d_y^k, d_z^k \right) \in \partial \mathcal{L}^\alpha(\mathbf{x}^k, y^k, z^k)$. By the norm property, we have $\|d_z^k\| = \frac{1}{\alpha\beta} \|\Delta z^k\|$, and by the triangle inequality from (12) and (13), we obtain

$$\|d_{x_i}^k\| \leq \|A_i\|\|\Delta z^k\| + \alpha\|A_i\|\Big(\sum_{j=i+1}^{p}\|A_j\|\|\Delta x_j^k\|\Big)$$
$$+\alpha\|A_i\|\|B\|\|\Delta y^k\| + \|Q_i\|\|\Delta x_i^k\|,$$
$$\|d_y^k\| \leq \|B\|\|\Delta z^k\| + \|P\|\|\Delta y^k\|.$$

Hence, we get

$$\||d^k\|| \leq \sum_{i=1}^{p}\|d_{x_i}^k\| + \|d_y^k\|$$
$$+\|d_z^k\| \leq \sum_{i=1}^{p}\|Q_i\|\|\Delta x_i^k\| + \sum_{i=2}^{p}\|A_i\|\|\Delta x_i^k\|\sum_{j=1}^{i-1}\|A_j\|$$
$$+\Big(\|P\| + \alpha\|B\|\sum_{i=1}^{p}\|A_i\|\Big)\|\Delta y^k\| + \Big(\sum_{i=1}^{p}\|A_i\| + \|B\| + \frac{1}{\alpha\beta}\Big)\|\Delta z^k\|$$
$$\leq \Big(\|\mathbf{Q}\| + \|\mathbf{A}\|^2\Big)\sum_{i=1}^{p}\|\Delta x_i^k\| + \Big(\|P\| + \alpha\|\mathbf{A}\|\|B\|\Big)\|\Delta y^k\|$$
$$+\Big(\|\mathbf{A}\| + \|B\| + \frac{1}{\alpha\beta}\Big)\|\Delta z^k\|$$
$$\leq \rho\Big(\sum_{i=1}^{p}\|\Delta x_i^k\| + \|\Delta y^k\| + \|\Delta z^k\|\Big).$$

□

**Lemma 3** (Limiting continuity) *Suppose that A1–A2 hold. If $(\mathbf{x}^*, y^*, z^*)$ is the limit point of $\{(\mathbf{x}^{k_j}, y^{k_j}, z^{k_j})\}_{j\geq 0}$, a subsequence of $\{\mathbf{x}^k, y^k, z^k\}_{k\geq 0}$ generated by the PADMM algorithm* (5). *Then $\mathcal{L}^\alpha(\mathbf{x}^*, y^*, z^*) = \lim_{j\to\infty}\mathcal{L}^\alpha(\mathbf{x}^{k_j}, y^{k_j}, z^{k_j})$.*

**Proof** Let $\{(\mathbf{x}^{k_j}, y^{k_j}, z^{k_j})\}_{j\geq 0}$ be a subsequence of the sequence $\{(\mathbf{x}^k, y^k, z^k)\}_{j\geq 0}$ generated by the PADMM algorithm such that $\lim_{j\to\infty}(\mathbf{x}^{k_j}, y^{k_j}, z^{k_j}) = (\mathbf{x}^*, y^*, z^*)$. By A1, for $i = 1, \ldots, p$ the functions $f_i$ are lower semicontinuous; hence,

$$f_i(x_i^*) \leq \liminf_{j\to\infty} f(x_i^{k_j}). \tag{15}$$

By the $x_i$-subproblem of the PADMM algorithm for any $v \in \mathbb{R}^{n_i}$, we have

$$\mathcal{L}^\alpha(x_{<i}^{k_j+1}, x_i^{k_j+1}, x_{>i}^{k_j}, y^{k_j}, z^{k_j}) + \frac{1}{2}\|\Delta x_i^{k_j+1}\|_{Q_i}^2$$
$$\leq \mathcal{L}^\alpha(x_{<i}^{k_j+1}, v, x_{>i}^{k_j}, y^{k_j}, z^{k_j}) + \frac{1}{2}\|v - x_i^{k_j}\|_{Q_i}^2.$$

We choose $v = x_i^*$ and exploit the cosine rule $\|a_1 + c\|^2 - \|a_2 + c\|^2 = \|a_2 - a_1\|^2 + 2\langle a_2 + c, a_2 - a_1\rangle$ with $a_1 = A_i x_i^*$, $a_2 = A_i x_i^{k_j+1}$ and $c = A_{<i} x_{<i}^{k_j+1} + A_{>i} x_{>i}^{k_j} +$

$By^{k_j} + b$ to obtain

$$f_i(x_i^{k_j+1}) - f_i(x_i^*) + \langle A_i^T z^{k_j}, x_i^{k_j+1} - x_i^* \rangle + \|\Delta x_i^{k_j+1}\|_{Q_i}^2 \leq \tfrac{\alpha}{2} \|A_i(x_i^{k_j+1} - x_i^*)\|^2$$
$$+ \alpha \Big\langle A_i(x_i^* - x_i^{k_j+1}), A_{<i} x_{<i}^{k_j+1} + A_i x_i^{k_j+1} + A_{>i} x_{>i}^{k_j} + By^{k_j} + b \Big\rangle.$$

Since $\{(\mathbf{x}^{k_j}, y^{k_j}, z^{k_j})\}_{j \geq 0}$ is a convergent sequence, it is bounded and the difference of consecutive terms approaches zero. Thus, the latter inequality reduces to $\limsup_{j \to \infty} f_i(x_i^{k_j+1}) - f_i(x_i^*) \leq 0$, and in view of (15), we then have $\lim_{j \to \infty} f_i(x_i^{k_j}) = f_i(x_i^*)$, which follows that $\lim_{j \to \infty} \sum_{i=1}^p f_i(x_i^{k_j}) = \sum_{i=1}^p f_i(x_i^*)$. By A2, $h(y)$ is smooth; hence, $\lim_{j \to \infty} h(y^{k_j}) = h(y^*)$. Clearly, $\mathbf{A}\mathbf{x}^{k_j} + By^{k_j} + b \to \mathbf{A}\mathbf{x}^* + By^* + b$ as $j \to \infty$. Thus, $\lim_{j \to \infty} \mathcal{L}^\alpha(\mathbf{x}^{k_j}, y^{k_j}, z^{k_j}) = \mathcal{L}^\alpha(\mathbf{x}^*, y^*, z^*)$. □

**Lemma 4** *(Limit point is critical point) The set of limit points of the sequence $\{(\mathbf{x}^k, y^k, z^k)\}_{k \geq 0}$ generated by the PADMM algorithm (5), denoted $\omega(\{(\mathbf{x}^k, y^k, z^k)\}_{k \geq 0})$, satisfies*

$$\omega(\{(\mathbf{x}^k, y^k, z^k)\}_{k \geq 0}) \subseteq \mathrm{crit}\,\mathcal{L}^\alpha$$
$$:= \left\{ (\mathbf{x}^*, y^*, z^*) : \begin{array}{l} 0 \in \partial f_i(x_i^*) + A_i^T z^*, \ i = 1, \ldots, p \\ \nabla h(y^*) + B^T z^* = 0, \\ \mathbf{A}\mathbf{x}^* + By^* + b = 0 \end{array} \right\},$$

*where $\mathrm{crit}\,\mathcal{L}^\alpha$ denotes the set of critical points of $\mathcal{L}^\alpha$.*

**Proof** Let $(\hat{\mathbf{x}}, \hat{y}, \hat{z}) \in \omega(\{(\mathbf{x}^k, y^k, z^k)\}_{k \geq 0})$, hence there exists a subsequence $\{(\mathbf{x}^{k_j}, y^{k_j}, z^{k_j})\}_{j \geq 0}$ of $\{(\mathbf{x}^k, y^k, z^k)\}_{k \geq 0}$ such that $\lim_{j \to \infty} (\mathbf{x}^{k_j}, y^{k_j}, z^{k_j}) = (\hat{\mathbf{x}}, \hat{y}, \hat{z})$. For the sequence $\{d^k\}_{k \geq 0}$, we have $d^{\bar{k}_j} \in \partial \mathcal{L}^\alpha(\mathbf{x}^{k_j}, y^{k_j}, z^{k_j})$ for $j \geq 1$ and $d^{k_j} \to 0$ as $j \to \infty$ while $(\mathbf{x}^{k_j}, y^{k_j}, z^{k_j}) \to (\hat{\mathbf{x}}, \hat{y}, \hat{z})$ and $\lim_{j \to \infty} \mathcal{L}^\alpha(\mathbf{x}^{k_j}, y^{k_j}, z^{k_j}) = \mathcal{L}^\alpha(\hat{\mathbf{x}}, \hat{y}, \hat{z})$. By the closeness criterion of the limiting subdifferential, we then have $0 \in \partial \mathcal{L}^\alpha(\hat{\mathbf{x}}, \hat{y}, \hat{z})$, and equivalently, $(\hat{\mathbf{x}}, \hat{y}, \hat{z}) \in \mathrm{crit}(\mathcal{L}^\alpha)$. □

**Lemma 5** *(bound dual) Suppose that A2 and A3 hold. Then for all $k \in \mathbb{Z}_+$ the following assertions are true.*

(i) $B^T z^k = \beta w^k + (1-\beta) B^T z^{k-1}$, *where* $w^k = -P \Delta y^k - \nabla h(y^k)$,
(ii) $\frac{1}{\alpha \beta} \|\Delta z^{k+1}\|^2 \leq c_4 \|\Delta y^{k+1}\|^2 + c_3 \|\Delta y^k\|^2 + c_5 \|B^T \Delta z^k\|^2 - c_5 \|B^T \Delta z^{k+1}\|^2$,
(iii) $\frac{1}{2\alpha} \|z^{k+1}\|^2 \leq c_1 \|P\|^2 \|\Delta y^{k+1}\|^2 + c_1 \|\nabla h(y^{k+1})\|^2 + c_6 \|B^T \Delta z^{k+1}\|^2$,

*where $c_i$, $i = 2, \ldots, 6$ are defined in (9).*

**Proof** (i) By the optimality condition of $y$-subproblem, we have

$$\alpha B^T (\mathbf{A}\mathbf{x}^k + By^k + b) + B^T z^{k-1} - w^k = 0$$

which together with the $z$-iterate $\mathbf{A}\mathbf{x}^k + By^k + b = \frac{1}{\alpha \beta}(z^k - z^{k-1})$ we obtain the result.

**(ii)** By A3 and the $z$ iterate of the PADMM algorithm, we have $\Delta z^{k+1} \in \text{Im}(B)$; thus, $\lambda_{++}^{B^T B} \|\Delta z^{k+1}\|^2 \leq \|B^T \Delta z^{k+1}\|^2$. Part (i) follows $B^T \Delta z^{k+1} = \beta \Delta w^{k+1} + (1-\beta) B^T \Delta z^k$ and equivalently

$$B^T \Delta z^{k+1} = \rho(\beta)\left(\frac{\beta \Delta w^{k+1}}{\rho(\beta)}\right) + |1-\beta|\left(\text{sign}(1-\beta) B^T \Delta z^k\right),$$

where $\text{sign}(\lambda) = 1$ if $\lambda \geq 0$ and $\text{sign}(\lambda) = -1$ if $\lambda < 0$. By $\rho(\beta) + |1-\beta| = 1$ and the convexity of $\|\cdot\|^2$, we have

$$\|B^T \Delta z^{k+1}\|^2 \leq \frac{\beta^2}{\rho(\beta)}\|\Delta w^{k+1}\|^2 + |1-\beta|\|B^T \Delta z^k\|^2. \tag{16}$$

This also follows that

$$\begin{aligned}
\rho(\beta)\lambda_{++}^{B^T B}\|\Delta z^{k+1}\|^2 &\leq \rho(\beta)\|B^T \Delta z^{k+1}\|^2 \\
&\leq \frac{\beta^2}{\rho(\beta)}\|\Delta w^{k+1}\|^2 \\
&\quad + |1-\beta|\|B^T \Delta z^k\|^2 - |1-\beta|\|B^T \Delta z^{k+1}\|^2.
\end{aligned}$$

Rearrange this to obtain

$$\begin{aligned}
\|\Delta z^{k+1}\|^2 &\leq \frac{\beta^2}{\rho(\beta)^2 \lambda_{++}^{B^T B}}\|\Delta w^{k+1}\|^2 + \frac{|1-\beta|}{\rho(\beta)\lambda_{++}^{B^T B}}\|B^T \Delta z^k\|^2 \\
&\quad - \frac{|1-\beta|}{\rho(\beta)\lambda_{++}^{B^T B}}\|B^T \Delta z^{k+1}\|^2.
\end{aligned} \tag{17}$$

By the optimality condition of $y$ subproblem, we have

$$\Delta w^{k+1} = -P\Delta y^{k+1} + P\Delta y^k + \nabla h(y^k) - \nabla h(y^{k+1}).$$

By $\|a+b\|^2 \leq 2\|a\|^2 + 2\|b\|^2$ and the fact that $\nabla h$ has $L_h$ Lipschitz gradients, we obtain

$$\|\Delta w^{k+1}\|^2 \leq 2(\|P\| + L_h)^2 \|\Delta y^{k+1}\|^2 + 2\|P\|^2 \|\Delta y^k\|^2.$$

Exploit this in (16) to get

$$\begin{aligned}
\|\Delta z^{k+1}\|^2 &\leq \frac{2\beta^2(\|P\| + L_h)^2}{\rho(\beta)^2 \lambda_{++}^{B^T B}}\|\Delta y^{k+1}\|^2 + \frac{2\beta^2 \|P\|^2}{\rho(\beta)^2 \lambda_{++}^{B^T B}}\|\Delta y^k\|^2 \\
&\quad + \frac{|1-\beta|}{\rho(\beta)\lambda_{++}^{B^T B}}\|B^T \Delta z^k\|^2 - \frac{|1-\beta|}{\rho(\beta)\lambda_{++}^{B^T B}}\|B^T \Delta z^{k+1}\|^2.
\end{aligned}$$

We divide both sides by $\alpha\beta$ to get

$$\frac{1}{\alpha\beta}\left\|\Delta z^{k+1}\right\|^2 \leq c_4\left\|\Delta y^{k+1}\right\|^2 + c_3\|\Delta y^k\|^2 + c_5\left\|B^T\Delta z^k\right\|^2 - c_5\left\|B^T\Delta z^{k+1}\right\|^2.$$

**(iii).** By part (i), we have $\beta B^T z^{k+1} = \beta w^{k+1} + (1-\beta)B^T(z^k - z^{k+1})$ which can be written equivalently as follows

$$\beta B^T z^{k+1} = \rho(\beta) \cdot \frac{\beta w^{k+1}}{\rho(\beta)} + |1-\beta| \cdot \text{sign}(1-\beta)B^T(z^k - z^{k+1}).$$

By $\rho(\beta) + |1-\beta| = 1$ and the convexity of $\|\cdot\|^2$, we get

$$\beta^2 \lambda_{++}^{B^T B} \|z^{k+1}\|^2 \leq \beta^2 \|B^T z^{k+1}\|^2 \leq \rho(\beta) \cdot \left\|\frac{\beta w^{k+1}}{\rho(\beta)}\right\|^2 + |1-\beta|\|B^T \Delta z^{k+1}\|^2.$$

Simplify and rearrange to get

$$\|z^{k+1}\|^2 \leq \frac{1}{\rho(\beta)\lambda_{++}^{B^T B}} \|w^{k+1}\|^2 + \frac{|1-\beta|}{\beta^2 \lambda_{++}^{B^T B}} \|B^T \Delta z^{k+1}\|^2.$$

Exploit $\|w^{k+1}\|^2 \leq 2\|P\|^2\|\Delta y^{k+1}\|^2 + 2\|\nabla h(y^{k+1})\|^2$ to obtain

$$\|z^{k+1}\|^2 \leq \frac{2\|P\|^2}{\rho(\beta)\lambda_{++}^{B^T B}} \|\Delta y^{k+1}\|^2 + \frac{2}{\rho(\beta)\lambda_{++}^{B^T B}} \|\nabla h(y^{k+1})\|^2$$
$$+ \frac{|1-\beta|}{\beta^2 \lambda_{++}^{B^T B}} \|B^T \Delta z^{k+1}\|^2.$$

Divide both sides by $2\alpha$ to obtain the result. This concludes the proof. □

**Lemma 6** *(change of $\mathcal{L}^\alpha$ during each variable update) Suppose that A1–A4 hold. The iterates in PADMM algorithm* (5) *satisfy*
- (i) $\mathcal{L}^\alpha(x_{<i}^{k+1}, x_i^k, x_{>i}^k, y^k, z^k) - \mathcal{L}^\alpha(x_{<i}^{k+1}, x_i^{k+1}, x_{>i}^k, y^k, z^k) \geq \frac{1}{2}\|\Delta x_i^{k+1}\|_{T_i}^2$ *where* $T_i = Q_i$ *or* $T_i = \alpha A_i^T A_i - \epsilon_i I$ *where* $\epsilon_i > 0$ *defined in A4 for* $i = 1, \ldots, p$.
- (ii) $\mathcal{L}^\alpha(\mathbf{x}^{k+1}, y^k, z^k) - \mathcal{L}^\alpha(\mathbf{x}^{k+1}, y^{k+1}, z^k) \geq \frac{1}{2}\|\Delta y^{k+1}\|_D^2$, *where D is defined in* (9).
- (iii) $\mathcal{L}^\alpha(\mathbf{x}^{k+1}, y^{k+1}, z^k) - \mathcal{L}^\alpha(\mathbf{x}^{k+1}, y^{k+1}, z^{k+1}) = -\frac{1}{\alpha\beta}\|\Delta z^{k+1}\|^2$.

**Proof** (i) Let $k \in \mathbb{Z}_+$ and $i \in \{1, \ldots, p\}$ be fixed. By the $x_i$ subproblem, we have

$$\mathcal{L}^\alpha(x_{<i}^{k+1}, x_i^k, x_{>i}^k, y^k, z^k) - \mathcal{L}^\alpha(x_{<i}^{k+1}, x_i^{k+1}, x_{>i}^k, y^k, z^k) \geq \|\Delta x_i^{k+1}\|_{Q_i}^2.$$

On the other hand, by the optimality condition of $x_i$-subminimization problem we have $\bar{d}_i \in \partial f_i(x_i^{k+1})$ where

$$\bar{d}_i := -A_i^T z^k - \alpha A_i^T \left(A_{<i} x_{<i}^{k+1} + A_i x_i^{k+1} + A_{>i} x_{>i}^k + B y^k + b\right).$$

By this, we then have

$$\mathcal{L}^{\alpha}(x_{<i}^{k+1}, x_i^k, x_{>i}^k, y^k, z^k) - \mathcal{L}^{\alpha}(x_{<i}^{k+1}, x_i^{k+1}, x_{>i}^k, y^k, z^k)$$
$$= f_i(x_i^k) - f_i(x_i^{k+1}) + \frac{\alpha}{2}\|A_i \Delta x_i^{k+1}\|^2 - \langle \bar{d}_i, \Delta x_i^{k+1}\rangle$$
$$\geq \frac{\alpha}{2}\|A_i \Delta x_i^{k+1}\|^2 - \frac{\epsilon_i}{2}\|\Delta x_i^{k+1}\|^2 = \frac{1}{2}\|\Delta x_i^{k+1}\|^2_{\alpha A_i^T A_i - \epsilon_i I}.$$

**(ii)** Let $k \geq 0$ be fixed. From the optimality condition of $y$ minimization problem, we have

$$\nabla h(y^{k+1}) + \alpha B^T(\mathbf{A}x^{k+1} + By^{k+1} + b) + B^T z^k + P\Delta y^{k+1} = 0.$$

Multiply this equation by $\Delta y^{k+1}$ and rearrange to obtain

$$-\alpha \langle B^T(\mathbf{A}x^{k+1} + By^{k+1} + b), \Delta y^{k+1}\rangle - \langle B^T z^k, \Delta y^{k+1}\rangle$$
$$= \langle \nabla h(y^{k+1}), \Delta y^{k+1}\rangle + \|\Delta y^{k+1}\|^2_P. \tag{18}$$

By this and the fact that $h(y)$ is $L_h$ Lipschitz continuous, we then get

$$\mathcal{L}^{\alpha}(x^{k+1}, y^k, z^k) - \mathcal{L}^{\alpha}(x^{k+1}, y^{k+1}, z^k)$$
$$= h(y^k) - h(y^{k+1}) + \frac{\alpha}{2}\|\mathbf{A}x^{k+1} + By^k + b\|^2 - \frac{\alpha}{2}\|\mathbf{A}x^{k+1} + By^{k+1}$$
$$+b\|^2 - \langle B^T z^k, \Delta y^{k+1}\rangle$$
$$= h(y^k) - h(y^{k+1}) + \frac{\alpha}{2}\|B\Delta y^{k+1}\|^2 - \alpha \langle \Delta y^{k+1}, B^T(\mathbf{A}x^{k+1} + By^{k+1} + b)\rangle$$
$$- \langle B^T z^k, \Delta y^{k+1}\rangle$$
$$= h(y^k) - h(y^{k+1}) + \langle \nabla h(y^{k+1}), \Delta y^{k+1}\rangle + \|\Delta y^{k+1}\|^2_P + \frac{\alpha}{2}\|B\Delta y^{k+1}\|^2$$
$$\geq \|\Delta y^{k+1}\|^2_P + \frac{\alpha}{2}\|B\Delta y^{k+1}\|^2 - \frac{L_h}{2}\|\Delta y^{k+1}\|^2 = \frac{1}{2}\|\Delta y^{k+1}\|^2_D.$$

**(iii)** By the $z$ iterate of the algorithm, we obtain

$$\mathcal{L}^{\alpha}(\mathbf{x}^{k+1}, y^{k+1}, z^k) - \mathcal{L}^{\alpha}(\mathbf{x}^{k+1}, y^{k+1}, z^{k+1})$$
$$= \langle z^k - z^{k+1}, \mathbf{A}x^{k+1} + By^{k+1} + b\rangle = -\frac{1}{\alpha\beta}\|\Delta z^{k+1}\|^2.$$

That completes the proof. □

**Lemma 7** (*change of $\mathcal{L}^{\alpha}$ after each iteration*) *We assume that A1–A4 hold. The iterates in PADMM algorithm* (5) *satisfy*

$$\mathcal{L}^{\alpha}(\mathbf{x}^{k+1}, y^{k+1}, z^{k+1}) + \frac{1}{2}\sum_{i=1}^{p}\|\Delta x_i^{k+1}\|^2_{T_i} + \frac{1}{2}\|\Delta y^{k+1}\|^2_{D-2\epsilon_0 c_4 I_m}$$

$$+ \frac{\epsilon_0 - 1}{\alpha \beta} \|\Delta z^{k+1}\|^2 + \epsilon_0 c_5 \|B^T \Delta z^{k+1}\|^2$$
$$\leq \mathcal{L}^\alpha(x^k, y^k, z^k) + \epsilon_0 c_5 \|B^T \Delta z^k\|^2 + \epsilon_0 c_3 \|\Delta y^k\|^2, \quad (19)$$

where $\epsilon_0 > 1$, and $c_3, c_4, c_5$ are defined in (9).

*Proof* By Lemma 6, we have

$$\mathcal{L}^\alpha(x^{k+1}, y^{k+1}, z^{k+1}) = \mathcal{L}^\alpha(x^{k+1}, y^{k+1}, z^k) + \frac{1}{\alpha \beta} \|\Delta z^{k+1}\|^2$$
$$\leq \mathcal{L}^\alpha(x^{k+1}, y^k, z^k) + \frac{1}{\alpha \beta} \|\Delta z^{k+1}\|^2 - \frac{1}{2} \|\Delta y^{k+1}\|_D^2$$
$$\leq \mathcal{L}^\alpha(x^k, y^k, z^k) + \frac{1}{\alpha \beta} \|\Delta z^{k+1}\|^2 - \frac{1}{2} \|\Delta y^{k+1}\|_D^2 - \frac{1}{2} \sum_{i=1}^{p} \|\Delta x_i^{k+1}\|_{T_i}^2.$$

By Lemma 5 for $\epsilon_0 > 1$, we have

$$\frac{\epsilon_0}{\alpha \beta} \|\Delta z^{k+1}\|^2$$
$$\leq \epsilon_0 c_4 \|\Delta y^{k+1}\|^2 + \epsilon_0 c_3 \|\Delta y^k\|^2 + \epsilon_0 c_5 \|B^T \Delta z^k\|^2 - \epsilon_0 c_5 \|B^T \Delta z^{k+1}\|^2.$$

Add this to the latter inequality and rearrange to obtain (19). □

We define a modified version of the augmented Lagrangian function (3). We consider $\bar{\mathcal{L}} : \mathbb{R}^n \times \mathbb{R}^q \times \mathbb{R}^m \times \mathbb{R}^q \times \mathbb{R}^m \to (-\infty, +\infty]$ defined by

$$\bar{\mathcal{L}}(x, y, z, y', z') = \mathcal{L}^\alpha(x, y, z) + \epsilon_0 c_5 \|B^T(z - z')\|^2 + \epsilon_0 c_3 \|y - y'\|^2, \quad (20)$$

where $\epsilon_0 > 1$, and $c_3$ and $c_5$ are given in (9). Evaluating $\bar{\mathcal{L}}$ at the point $(\mathbf{x}^k, y^k, z^k, y^{k-1}, z^{k-1})$ where $\{(\mathbf{x}^k, y^k, z^k)\}_{k \geq 1}$ is the sequence generated by the PADMM algorithm (5) gives

$$\bar{\mathcal{L}}_k := \bar{\mathcal{L}}(\mathbf{x}^k, y^k, z^k, y^{k-1}, z^{k-1}) = \mathcal{L}^\alpha(\mathbf{x}^k, y^k, z^k) + \epsilon_0 c_5 \|B^T \Delta z^k\|^2 + \epsilon_0 c_3 \|\Delta y^k\|^2. \quad (21)$$

**Lemma 8** *Suppose that A1–A4 hold, and*

$$\bar{\tau} = \min\{\max\{q_i, \lambda_{++}^{\alpha A_i^T A_i - \epsilon_i I}\} : i = 1, \ldots, p\}.$$

*If $\sigma = \min\{\lambda_{\min}(\bar{D}), \bar{\tau}, \frac{\epsilon_0 - 1}{\alpha \beta}\} > 0$ where $\bar{D}$ is defined in (9), then the following assertions hold.*

*(i)* The sequence $\{\bar{\mathcal{L}}_k\}_{k \geq 1}$ is monotonically decreasing and

$$\bar{\mathcal{L}}_{k+1} + \sigma \Big[ \sum_{i=1}^{p} \|\Delta x_i^{k+1}\|^2 + \|\Delta y^{k+1}\|^2 + \|\Delta z^{k+1}\|^2 \Big] \leq \bar{\mathcal{L}}_k. \quad (22)$$

*(ii)* If the sequence $\{\mathbf{x}^k, y^k, z^k\}_{k \geq 0}$ generated by the PADMM is bounded then $\bar{\mathcal{L}}_k$ is bounded from below for all $k \in \mathbb{N}$ and converges as k approaches infinity.

**Proof** **(i)** By (19) and (21), we obtain

$$\bar{\mathcal{L}}_{k+1} + \frac{1}{2} \sum_{i=1}^{p} \|\Delta x_i^{k+1}\|_{T_i}^2 + \frac{1}{2} \|\Delta y^{k+1}\|_{D - 2\epsilon_0 (c_3 + c_4)}^2 + \frac{\epsilon_0 - 1}{\alpha \beta} \|\Delta z^{k+1}\|^2 \leq \bar{\mathcal{L}}_k.$$

We have $T_i = Q_i \succeq q_i I$ where $q_i > 0$ or $T_i = \alpha A_i^T A_i - \epsilon_i I$. We let $\bar{\tau} = \min\{\max\{q_i, \lambda_{++}^{\alpha A_i^T A_i - \epsilon_i I}\} : i = 1, \ldots, p\}$ and $\bar{D} = D - 2\epsilon_0(c_3 + c_4)I$. For $\sigma = \min\{\lambda_{\min}(\bar{D}), \bar{\tau}, \frac{\epsilon_0 - 1}{\alpha \beta}\} > 0$ we arrive at the inequality (22), showing monotonicity of $\{\bar{\mathcal{L}}_k\}_{k \geq 1}$.

**(ii)** If $\{\mathbf{x}^k, y^k, z^k\}_{k \geq 0}$ is a bounded sequence, then $\|B^T \Delta z^k\|^2 + \|\Delta y^k\|^2$ is also bounded, and bounded from below by zero. Hence, we only need to show that $\mathcal{L}^\alpha(\mathbf{x}^k, y^k, z^k)$ is bounded. By the boundedness of $\{(\mathbf{x}^k, y^k, z^k)\}_{k \geq 0}$, the terms $\langle z^k, \mathbf{A}\mathbf{x}^k + By^k + b \rangle$ and $\|\mathbf{A}\mathbf{x}^k + By^k + b\|^2$ are bounded. By A1 the functions $f_i, i = 1, \ldots, p$ are coercive, since $\{\mathbf{x}^k\}_{k \geq 0}$ is bounded, then $\{f_i(x_i^k)\}_{k \geq 0}$ is bounded for $i = 1, \ldots, p$. By A2, the function $h$ is bounded from below. Thus, $\{\mathcal{L}^\alpha(\mathbf{x}^k, y^k, z^k)\}_{k \geq 0}$ is bounded from below. Since $\{\bar{\mathcal{L}}_k\}_{k \geq 1}$ is bounded from below and by part (i) it is monotonically decreasing then we conclude that $\{\bar{\mathcal{L}}_k\}_{k \geq 1}$ converges. □

**Lemma 9** *Suppose that A1–A4 hold and $(\mathbf{x}^*, y^*, z^*) \in \omega(\{(\mathbf{x}^k, y^k, z^k)\}_{k \geq 0})$, then*

$$\bar{\mathcal{L}}(\mathbf{x}^*, y^*, z^*, y^*, z^*) = \mathcal{L}^\alpha(\mathbf{x}^*, y^*, z^*) = \sum_{i=1}^{p} f_i(x_i^*) + h(y^*).$$

**Proof** Let $(\mathbf{x}^*, y^*, z^*) \in \omega(\{(\mathbf{x}^k, y^k, z^k)\}_{k \geq 0})$ and let the subsequence $\{(\mathbf{x}^{k_j}, y^{k_j}, z^{k_j})\}_{j \geq 0}$ of the PADMM algorithm converges to $(\mathbf{x}^*, y^*, z^*)$. Clearly $\|\Delta y^{k_j}\|$ and $\|B^T \Delta z^{k_j}\|$ both converge to zero as $j$ goes to infinity. This, together with Lemmas 3 and 4, gives

$$\lim_{j \to \infty} \bar{\mathcal{L}}_{k_j} = \lim_{j \to \infty} \mathcal{L}^\alpha(x^{k_j}, y^{k_j}, z^{k_j}) = \mathcal{L}^\alpha(x^*, y^*, z^*)$$

$$= \sum_{i=1}^{p} f_i(x_i^*) + h(y^*) = \bar{\mathcal{L}}(\mathbf{x}^*, y^*, z^*, y^*, z^*).$$

□

**Lemma 10** *(subgradient bound for $\bar{\mathcal{L}}_k$) Let $\{(\mathbf{x}^k, y^k, z^k)\}_{k \geq 0}$ be a sequence generated by the PADMM algorithm and assume that $(\{d_{x_i}^k\}_{i=1}^p, d_y^k, d_z^{\bar{k}}) \in \partial \mathcal{L}^\alpha(\mathbf{x}^k, y^k, z^k)$. Then*

$$s^k := \left(\{s_{x_i}^k\}_{i=1}^p, s_y^k, s_z^k, s_{y'}^k, s_{z'}^k\right) \in \partial \bar{\mathcal{L}}(\mathbf{x}^k, y^k, z^k, y^{k-1}, z^{k-1}), \quad k \geq 1$$

*where $s_{x_i}^k := d_{x_i}^k$, $i = 1, \ldots, p$, $s_y^k := d_y^k + 2\epsilon_0 c_3 \Delta y^k$, $s_z^k := d_z^k + 2\epsilon_0 c_5 B B^T \Delta z^k$, $s_{y'}^k := -2\epsilon_0 c_3 \Delta y^k$, $s_{z'}^k := -2\epsilon_0 c_5 B B^T \Delta z^k$, and $\epsilon_0 > 1$, and $c_3$ and $c_5$ are defined in (9). Moreover,*

$$|||s^k||| \leq \tilde{\rho}\left(\sum_{i=1}^p \|\Delta x_i^k\| + \|\Delta y^k\| + \|\Delta z^k\|\right), \tag{23}$$

*where $\tilde{\rho} := \max\left\{\|\mathbf{Q}\| + \|\mathbf{A}\|^2, \|P\| + \alpha\|\mathbf{A}\|\|B\| + 4\epsilon_0 c_3, \|\mathbf{A}\| + \|B\| + 4\epsilon_0 c_5 \|B\|^2 + \frac{1}{\alpha\beta}\right\}$, and $\|\mathbf{Q}$ and $\|\mathbf{A}\|$ are defined in Lemma 2.*

**Proof** Let $k \geq 1$ be fixed and $(\{d_{x_i}^k\}_{i=1}^p, d_y^k, d_z^k) \in \partial \mathcal{L}^\alpha(x^k, y^k, z^k)$. By taking partial derivatives of $\bar{\mathcal{L}}$ with respect to $x, y, z, y', z'$ and evaluating them at $(\mathbf{x}^k, y^k, z^k, y^{k-1}, z^{k-1})$, we obtain

$$\partial_{x_i} \bar{\mathcal{L}}(\mathbf{x}^k, y^k, z^k, y^{k-1}, z^{k-1}) = \partial_{x_i} \mathcal{L}^\alpha(\mathbf{x}^k, y^k, z^k) \ni d_{x_i}^k, \quad i = 1, \ldots, p;$$

$$\nabla_y \bar{\mathcal{L}}(\mathbf{x}^k, y^k, z^k, y^{k-1}, z^{k-1}) = \nabla_y \mathcal{L}^\alpha(\mathbf{x}^k, y^k, z^k) + 2\epsilon_0 c_3 \Delta y^k = d_y^k + 2\epsilon_0 c_3 \Delta y^k;$$

$$\nabla_z \bar{\mathcal{L}}(\mathbf{x}^k, y^k, z^k, y^{k-1}, z^{k-1}) = \nabla_z \mathcal{L}^\alpha(\mathbf{x}^k, y^k, z^k)$$
$$\qquad\qquad\qquad + 2\epsilon_0 c_5 B B^T \Delta z^k = d_z^k + 2\epsilon_0 c_5 B B^T \Delta z^k;$$

$$\nabla_{y'} \bar{\mathcal{L}}(\mathbf{x}^k, y^k, z^k, y^{k-1}, z^{k-1}) = -2\epsilon_0 c_3 \Delta y^k;$$

$$\nabla_{z'} \bar{\mathcal{L}}(\mathbf{x}^k, y^k, z^k, y^{k-1}, z^{k-1}) = -2\epsilon_0 c_5 B B^T \Delta z^k.$$

These equations follow that

$$s^k = \begin{bmatrix} s_{x_1}^k \\ \vdots \\ s_{x_p}^k \\ s_y^k \\ s_z^k \\ s_{y'}^k \\ s_{z'}^k \end{bmatrix} = \begin{bmatrix} d_{x_1}^k \\ \vdots \\ d_{x_p}^k \\ d_y^k + 2\epsilon_0 c_3 \Delta y^k \\ d_z^k + 2\epsilon_0 c_5 B B^T \Delta z^k \\ -2\epsilon_0 c_3 \Delta y^k \\ -2\epsilon_0 c_5 B B^T \Delta z^k \end{bmatrix} \in \partial \bar{\mathcal{L}}(\mathbf{x}^k, y^k, z^k, y^{k-1}, z^{k-1}).$$

By the triangle inequality, we have

$$\|s_y^k\| \leq \|d_y^k\| + 2\epsilon_0 c_3 \|\Delta y^k\|, \qquad \|s_z^k\| \leq \|d_z^k\| + 2\epsilon_0 c_5 \|B\|^2 \|\Delta z^k\|,$$

$$\|s_{y'}^k\| \leq 2\epsilon_0 c_3 \|\Delta y^k\|, \qquad \|s_z^k\| \leq 2\epsilon_0 c_5 \|B\|^2 \|\Delta z^k\|.$$

Thus,

$$|||s^k||| \leq \sum_{i=1}^{p} \|s_{x_i}^k\| + \|s_y^k\| + \|s_z^k\| + \|s_{y'}^k\| + \|s_z^k\|$$

$$\leq \sum_{i=1}^{p} \|d_{x_i}^k\| + \|d_y^k\| + \|d_z^k\| + 4\epsilon_0 c_3 \|\Delta y^k\| + 4\epsilon_0 c_5 \|B\|^2 \|\Delta z^k\|.$$

By Lemma 2, we have

$$|||s^k||| \leq \left(\|\mathbf{Q}\| + \|\mathbf{A}\|^2\right) \sum_{i=1}^{p} \|\Delta x_i^k\| + \left(\|P\| + \alpha\|\mathbf{A}\|\|B\| + 4\epsilon_0 c_3\right) \|\Delta y^k\|$$

$$+ \left(\|\mathbf{A}\| + \|B\| + 4\epsilon_0 c_5 \|B\|^2 + \frac{1}{\alpha\beta}\right) \|\Delta z^k\| \leq \tilde{\rho}\Big(\sum_{i=1}^{p} \|\Delta x_i^k\| + \|\Delta y^k\| + \|\Delta z^k\|\Big).$$

□

## 4 Convergence

In this section, we prove that the sequence generated with PADMM algorithm (5) converges. We begin with some useful lemmas.

**Lemma 11** *Suppose that A1–A4 hold and $\sigma$ defined in Lemma 8 is positive. Then the sequence $\{\mathbf{x}^k, y^k, z^k\}_{k\geq 0}$ generated by the PADMM algorithm (5) is bounded, and $\lim_{k\to\infty} \|\Delta x_1^k\| = \cdots = \lim_{k\to\infty} \|\Delta x_p^k\| = \lim_{k\to\infty} \|\Delta y^k\| = \lim_{k\to\infty} \|\Delta z^k\| = 0$.*

***Proof*** Let $\{(\mathbf{x}^k, y^k, z^k)\}_{k\geq 0}$ be a sequence generated by the PADMM algorithm. By $\sigma > 0$, the sequence $\{\bar{\mathcal{L}}_k\}_{k\geq 1}$ defined in (21) is monotonically decreasing; hence, it is bounded above by $\bar{\mathcal{L}}_1 = \bar{\mathcal{L}}(\mathbf{x}^1, y^1, z^1, y^0, z^0)$. This follows that

$$\sum_{i=1}^{p} f(x_i^k) + h(y^k) + \frac{\alpha}{2}\|\mathbf{A}\mathbf{x}^k + By^k + \alpha^{-1}z^k + b\|^2 - \frac{1}{2\alpha}\|z^k\|^2$$
$$+ \sigma \sum_{i=1}^{p} \|\Delta x_i^k\|^2 + \sigma \|\Delta z^k\|^2 + (\sigma + \epsilon_0 c_3)\|\Delta y^k\|^2 + \epsilon_0 c_5 \|B^T \Delta z^k\|^2 \leq \bar{\mathcal{L}}_1 \quad (24)$$

for all $k \in \mathbb{Z}_+$. By Lemma 5(iii), (24) leads to

$$\sum_{i=1}^{p} f(x_i^k) + \frac{\alpha}{2}\|\mathbf{A}\mathbf{x}^k + By^k + \alpha^{-1}z^k + b\|^2 + r_1\|\Delta y^k\|^2$$
$$+ \sigma \sum_{i=1}^{p} \|\Delta x^k\|^2 + r_2 \|B^T \Delta z^k\|^2 + \sigma\|\Delta z^k\|^2 \leq \bar{\mathcal{L}}_1 - \inf_y \left\{h(y) - c_1\|\nabla h(y)\|^2\right\}, \quad (25)$$

where $r_1 := \epsilon_0 c_3 + \sigma - \|P\|^2 c_1 > 0$ and $r_2 := \epsilon_0 c_5 - c_6 > 0$. By A2, since $h$ is $L_h$ Lipschitz differentiable then for any $k \geq 1$ and $\delta > 0$ we have $h(y^k - \delta \nabla h(y^k)) \leq h(y^k) - (\delta - \frac{L_h \delta^2}{2})\|\nabla h(y^k)\|^2$. Again by A2, $h$ is bounded from below then

$$-\infty < \inf\{h(y) - \left(\delta - \frac{L_h \delta^2}{2}\right)\|\nabla h(y)\|^2 \,:\, y \in \mathbb{R}^q\}. \tag{26}$$

We choose $\delta(> 0)$ such that $\delta - \frac{L_h \delta^2}{2} = c_1$. By (26), the right-hand side of (25) is finite; hence,

$$\sum_{i=1}^{p} f(x_i^k) + \|\mathbf{A}\mathbf{x}^k + By^k + \alpha^{-1}z^k + b\|^2 + \|\Delta y^k\|^2 + \|\Delta z^k\|^2 < +\infty. \tag{27}$$

By A1, since $f_i$, $i = 1, \ldots, p$ are coercive and $\sum_{i=1}^{p} f(x_i^k) < +\infty$, then the sequence $\{x_i^k\}_{k \geq 0}$ is bounded for $i = 1, \ldots, p$, and so are $\{\mathbf{x}^k\}_{k \geq 0}$ and $\{\mathbf{A}\mathbf{x}^k\}_{k \geq 0}$.

By the $z$ iterate of PADMM we have $By^k = \frac{1}{\alpha \beta} \Delta z^k - \mathbf{A}\mathbf{x}^k - b$. By (27), $\{\Delta z^k\}_{k \geq 0}$ is bounded. By A3, $B^T B$ is invertible; thus, $\{y^k\}_{k \geq 0}$ is bounded. Finally, by (27) $\{\mathbf{A}\mathbf{x}^k + By^k + \frac{1}{\alpha} z^k\}_{k \geq 0}$ is bounded and $\mathbf{A}\mathbf{x}^k$ and $By^k$ are bounded, so $\{z^k\}_{k \geq 0}$ is then bounded. Summing up (22) for $k = 1, \ldots, K-1$, $K \geq 2$ gives

$$\sum_{k=1}^{K} \left\{ \sum_{i=1}^{p} \|\Delta x_i^k\|^2 + \|\Delta y^k\|^2 + \|\Delta z^k\|^2 \right\} \leq \frac{1}{\sigma}(\bar{\mathcal{L}}_1 - \bar{\mathcal{L}}_K).$$

Since $\{(\mathbf{x}^k, y^k, z^k)\}_{k \geq 0}$ is bounded, by Lemma 8 (ii) the right-hand side is finite for any $K$. We let $K$ approach to infinity to obtain

$$\sum_{k=1}^{\infty} \left\{ \sum_{i=1}^{p} \|\Delta x_i^k\|^2 + \|\Delta y^k\|^2 + \|\Delta z^k\|^2 \right\} < \frac{1}{\sigma}(\bar{\mathcal{L}}_1 - \inf_{k \geq 1} \bar{\mathcal{L}}_k) < +\infty.$$

This follows $\lim_{k \to \infty} \left( \sum_{i=1}^{p} \|\Delta x_i^k\|^2 + \|\Delta y^k\|^2 + \|\Delta z^k\|^2 \right) = 0$, and consequently $\lim_{k \to \infty} \|\Delta x_1^k\| = \cdots = \lim_{k \to \infty} \|\Delta x_p^k\| = \lim_{k \to \infty} \|\Delta y^k\| = \lim_{k \to \infty} \|\Delta z^k\| = 0$. □

**Lemma 12** *Let $\{(\mathbf{x}^k, y^k, z^k)\}_{k \geq 0}$ be a sequence generated by the PADMM algorithm (5), and denote*

$$\Omega := \omega\left(\{(\mathbf{x}^k, y^k, z^k, y^{k-1}, z^{k-1})\}_{k \geq 1}\right). \tag{28}$$

*Then the following statements are true*

(i) $\Omega$ *is nonempty, connected, and compact.*
(ii) $\Omega \subseteq \{(\mathbf{x}^*, y^*, z^*, y^*, z^*) \in \mathbb{R}^n \times \mathbb{R}^q \times \mathbb{R}^m \times \mathbb{R}^q \times \mathbb{R}^m \,:\, (\mathbf{x}^*, y^*, z^*) \in \mathrm{crit}(\mathcal{L}^\alpha)\}$

(iii) $\lim_{k\to\infty} \text{dist}\big[(\mathbf{x}^k, y^k, z^k, y^{k-1}, x^{k-1}), \Omega\big] = 0$.
(iv) The sequences $\{\bar{\mathcal{L}}_k\}_{k\geq 1}$, $\{\bar{\mathcal{L}}^\alpha(\mathbf{x}^k, y^k, z^k)\}_{k\geq 0}$, and $\{\sum_{i=1}^p f_i(x_i^k) + h(y^k)\}_{k\geq 0}$ approaches the same limit on $\Omega$.

*Proof* These results follow from Lemmas 8, 9, and 11. Thus, we omit the proof. □

**Theorem 1** *Suppose that A1–A4 hold and $\sigma$ defined in Lemma 8 is positive. Let the sequence $\{(\mathbf{x}^k, y^k, z^k)\}_{k\geq 0}$ generated by the PADMM algorithm (5) be bounded and $\bar{\mathcal{L}}(\mathbf{x}, y, z, y', z')$ satisfies the KŁ property on $\Omega$ (28); that is, for every $v^* := (\mathbf{x}^*, y^*, z^*, y^*, z^*) \in \Omega$ there exists $\epsilon > 0$, $\eta \in [0, +\infty)$, and a function $\psi \in \Psi_\eta$ such that for any $v := (\mathbf{x}, y, z, y', z') \in \mathcal{S}$ where*

$$\mathcal{S} := \Big\{ v \in \mathbb{R}^n \times \mathbb{R}^m \times \mathbb{R}^p \times \mathbb{R}^m \times \mathbb{R}^p : \\ \text{dist}(v, \Omega) < \epsilon \Big\} \cap [\bar{\mathcal{L}}(v^*) < \bar{\mathcal{L}} < \bar{\mathcal{L}}(v^*) + \eta], \quad (29)$$

*it holds*

$$\psi'\Big(\bar{\mathcal{L}}(v) - \bar{\mathcal{L}}(v^*)\Big) \cdot \text{dist}\big(0, \partial\bar{\mathcal{L}}(v)\big) \geq 1. \quad (30)$$

*Then $\{u^k\}_{k\geq 0} := \{(\mathbf{x}^k, y^k, z^k)\}_{k\geq 0}$ satisfies*

$$\sum_{k=1}^\infty |||\Delta u^k||| \leq \sum_{k=1}^\infty \Big\{ \sum_{i=1}^p \|\Delta x_i^k\| + \|\Delta y^k\| + \|\Delta z^k\| \Big\} < +\infty,$$

*and consequently it converges to a stationary point of (2).*

*Proof* By Lemma 8, the sequence $\{\bar{\mathcal{L}}_k\}_{k\geq 1}$ is monotonically decreasing and converges. Let $\bar{\mathcal{L}}_* := \lim_{k\to\infty} \bar{\mathcal{L}}_k$, then the error sequence $e_k := \bar{\mathcal{L}}_k - \bar{\mathcal{L}}_*$, $k \geq 1$ is nonnegative, monotonically decreasing, and converges to zero. We prove under two cases.

*Case 1:* let assume that there is a $k_1 \in \mathbb{Z}_+$ such that $e_{k_1} = 0$. Since $\{e_k\}_{k\geq 1}$ is monotonically decreasing, $e_k = 0$ for all $k \geq k_1$. By Lemma 8(i), then we have

$$\sum_{i=1}^p \|\Delta x_i^{k+1}\|^2 + \|\Delta y^{k+1}\|^2 + \|\Delta z^{k+1}\|^2 \leq \frac{1}{\sigma}(e_k - e_{k+1}) = 0, \quad \forall k \geq k_1.$$

Hence, as $\{(\mathbf{x}^k, y^k, z^k)\}_{k\geq 0}$ is a bounded sequence then

$$\sum_{k=1}^\infty \Big\{ \sum_{i=1}^p \|\Delta x_i^k\| + \|\Delta y^k\| + \|\Delta z^k\| \Big\} = \sum_{k=1}^{k_1} \Big\{ \sum_{i=1}^p \|\Delta x_i^k\| + \|\Delta y^k\| + \|\Delta z^k\| \Big\}$$

is bounded.

*Case 2:* we let $e_k > 0$ for all $k \geq 1$. By Lemma 8(i) we have

$$\sum_{i=1}^{p} \|\Delta x_i^{k+1}\|^2 + \|\Delta y^{k+1}\|^2 + \|\Delta z^{k+1}\|^2 \leq \frac{1}{\sigma}(e_k - e_{k+1}), \quad \forall k \geq 1. \quad (31)$$

By Lemma 12, $\Omega$ is nonempty, compact, and connected, and $\bar{\mathcal{L}}$ takes on a constant value $\bar{\mathcal{L}}_*$ on $\Omega$. Since $\bar{\mathcal{L}}_k \searrow \bar{\mathcal{L}}_*$ as $k \to \infty$ and $\bar{\mathcal{L}}$ satisfies the KŁ property, then there exists $k_1 \geq 1$ such that $\bar{\mathcal{L}}_* < \bar{\mathcal{L}}_k < \bar{\mathcal{L}}_* + \eta$ for all $k \geq k_1$ and some $\eta > 0$. Moreover, $\lim_{k\to\infty} \text{dist}\big[(\mathbf{x}^k, y^k, z^k, y^{k-1}, x^{k-1}), \Omega\big] = 0$ by Lemma 12; thus, for a given $\epsilon > 0$ there exists $k_2 \geq 1$ such that for all $k \geq k_2$ it holds

$$\text{dist}\Big[(\mathbf{x}^k, y^k, z^k, y^{k-1}, x^{k-1}), \Omega\Big] < \epsilon.$$

Let $\tilde{k} := \max\{k_1, k_2, 3\}$, then $(\mathbf{x}^k, y^k, z^k, y^{k-1}, x^{k-1}) \in \mathcal{S}$ and $\psi'(\bar{\mathcal{L}}_k - \bar{\mathcal{L}}_*) \cdot \text{dist}\big(0, \partial \bar{\mathcal{L}}_k\big) \geq 1$ for all $k \geq \tilde{k}$, or equivalently

$$\psi'(e_k) \cdot \text{dist}\Big(0, \partial \bar{\mathcal{L}}_k\Big) \geq 1. \quad (32)$$

Multiply both sides of this inequality by $|||\Delta u^{k+1}|||^2$ to get

$$\psi'(e_k) \cdot \text{dist}\Big(0, \partial \bar{\mathcal{L}}_k\Big)|||\Delta u^{k+1}|||^2 \geq |||\Delta u^{k+1}|||^2, \quad (33)$$

where $|||\Delta u^k|||^2 = \sum_{i=1}^{p} \|\Delta x_i^k\|^2 + \|\Delta y^k\|^2 + \|\Delta z^k\|^2$. By (31), (33) leads to

$$\frac{1}{\sigma}\psi'(e_k)(e_k - e_{k+1}) \cdot \text{dist}\Big(0, \partial \bar{\mathcal{L}}_k\Big) \geq |||\Delta u^{k+1}|||^2. \quad (34)$$

By the concavity of $\psi$, $\psi(e_k) - \psi(e_{k+1}) \geq \psi'(e_k)(e_k - e_{k+1})$; hence,

$$|||\Delta u^{k+1}|||^2 \leq \frac{1}{\sigma}\big(\psi(e_k) - \psi(e_{k+1})\big) \cdot \text{dist}(0, \partial \bar{\mathcal{L}}_k). \quad (35)$$

By the arithmetic mean-geometric mean inequality, for any $\gamma > 0$ we have

$$\sqrt{\frac{1}{\sigma}\big(\psi(e_k) - \psi(e_{k+1})\big) \cdot \text{dist}(0, \partial \bar{\mathcal{L}}_k)} \leq \frac{\gamma}{2\sigma}\big(\psi(e_k) - \psi(e_{k+1})\big) + \frac{1}{2\gamma}\text{dist}(0, \partial \bar{\mathcal{L}}_k). \quad (36)$$

Moreover,

$$\frac{1}{\sqrt{p+2}}\Big(\sum_{i=1}^{p} \|\Delta x_i^{k+1}\| + \|\Delta y^{k+1}\| + \|\Delta z^{k+1}\|\Big) \leq |||\Delta u^{k+1}|||. \quad (37)$$

Thus, (35), (36), and (37) give rise to

$$\sum_{i=1}^{p} \|\Delta x_i^{k+1}\| + \|\Delta y^{k+1}\| + \|\Delta z^{k+1}\|$$
$$\leq \frac{\gamma\sqrt{p+2}}{2\sigma}\bigl(\psi(e_k) - \psi(e_{k+1})\bigr) + \frac{\sqrt{p+2}}{2\gamma}\text{dist}(0, \partial\bar{\mathcal{L}}_k). \quad (38)$$

Note that $\text{dist}(0, \partial\bar{\mathcal{L}}_k) = |||s^k|||$, where $s^k \in \partial\bar{\mathcal{L}}_k$ defined in Lemma 10, hence

$$\sum_{i=1}^{p} \|\Delta x_i^{k+1}\| + \|\Delta y^{k+1}\| + \|\Delta z^{k+1}\|$$
$$\leq \frac{\gamma\sqrt{p+2}}{2\sigma}\bigl(\psi(e_k) - \psi(e_{k+1})\bigr) + \frac{\tilde{\rho}\sqrt{p+2}}{2\gamma}\Bigl(\sum_{i=1}^{p} \|\Delta x_i^k\| + \|\Delta y^k\| + \|\Delta z^k\|\Bigr). \quad (39)$$

For any $K \geq k \geq \underline{k}$, the following identity holds

$$\sum_{k=\underline{k}}^{K} \|\Delta w^k\| = \sum_{k=\underline{k}}^{K} \|\Delta w^{k+1}\| + \|\Delta w^{\underline{k}}\| - \|\Delta w^{K+1}\|.$$

Choose $\gamma > 0$ large enough such that $2\gamma > \tilde{\rho}\sqrt{p+2}$. Let $\delta_0 = 1 - \frac{\tilde{\rho}\sqrt{p+2}}{2\gamma} > 0$. Summing up (39) from $k = \underline{k} \geq \tilde{k}$ to $K \geq \underline{k}$ leads to

$$\sum_{k=\underline{k}}^{K} \Bigl\{\sum_{i=1}^{p} \|\Delta x_i^{k+1}\| + \|\Delta y^{k+1}\| + \|\Delta z^{k+1}\|\Bigr\} \leq \delta_1\bigl(\psi(e_{\underline{k}}) - \psi(e_{K+1})\bigr)$$
$$+ \delta_2\Bigl(\sum_{i=1}^{p} \|\Delta x_i^{\underline{k}}\| + \|\Delta y^{\underline{k}}\| + \|\Delta z^{\underline{k}}\|\Bigr)$$
$$- \delta_2\Bigl(\sum_{i=1}^{p} \|\Delta x_i^{K+1}\| + \|\Delta y^{K+1}\| + \|\Delta z^{K+1}\|\Bigr),$$

where $\delta_1 := \gamma\sqrt{p+2}/2\sigma\delta_0$ and $\delta_2 = \tilde{\rho}\sqrt{p+2}/2\gamma\delta_0$. Recall that $e_k$ is monotonically decreasing and $\psi' > 0$, therefore

$$\sum_{k=\underline{k}}^{K} \Bigl\{\sum_{i=1}^{p} \|\Delta x_i^{k+1}\| + \|\Delta y^{k+1}\| + \|\Delta z^{k+1}\|\Bigr\}$$
$$\leq \delta_1\psi(e_{\underline{k}}) + \delta_2\Bigl(\sum_{i=1}^{p} \|\Delta x_i^{\underline{k}}\| + \|\Delta y^{\underline{k}}\| + \|\Delta z^{\underline{k}}\|\Bigr).$$

The right-hand side of this inequality is bounded for any $K \geq \underline{k}$, we let $K \to \infty$ to obtain

$$\sum_{k=\underline{k}}^{\infty} \left( \sum_{i=1}^{p} \|\Delta x_i^{k+1}\| + \|\Delta y^{k+1}\| + \|\Delta z^{k+1}\| \right)$$
$$\leq \delta_1 \psi(e_{\underline{k}}) + \delta_2 \left( \sum_{i=1}^{p} \|\Delta x_i^{\underline{k}}\| + \|\Delta y^{\underline{k}}\| + \|\Delta z^{\underline{k}}\| \right). \quad (40)$$

Then (40) together with the fact that $\{(\mathbf{x}^k, y^k, z^k)\}_{k \geq 0}$ is bounded gives

$$\sum_{k=0}^{\infty} \left\{ \sum_{i=1}^{p} \|\Delta x_i^{k+1}\| + \|\Delta y^{k+1}\| + \|\Delta z^{k+1}\| \right\} < +\infty. \quad (41)$$

For any $p, q, K \in \mathbb{Z}_+$ where $q \geq p$, we have

$$|||u^q - u^p||| = |||\sum_{k=p}^{q-1} \Delta u^{k+1}||| \leq \sum_{k=p}^{q-1} |||\Delta u^{k+1}||| \leq \sum_{k=1}^{\infty} |||\Delta u^{k+1}|||.$$

Since by (41) the right-hand side is finite, then

$$\sum_{k=1}^{\infty} |||\Delta u^{k+1}||| \leq \sum_{k=0}^{\infty} \left\{ \sum_{i=1}^{p} \|\Delta x_i^{k+1}\| + \|\Delta y^{k+1}\| + \|\Delta z^{k+1}\| \right\} < +\infty.$$

This implies that $\{u^k\}_{k \geq 0} = \{(\mathbf{x}, y^k, z^k)\}_{k \geq 0}$ is a Cauchy sequence. By Lemma 3, it converges to a stationary point. □

**Remark 1** If $Q_i, i = 1, \ldots, p$ and $P$ are chosen zero matrices, the PADMM algorithm (5) is reduced to the ADMM algorithm (4). If $\alpha > \min\{\frac{\epsilon_1}{\|A_1\|^2}, \ldots, \frac{\epsilon_p}{\|A_p\|^2}, \frac{L_h}{\|B\|^2}\}$ where $\epsilon_i > 0, i = 1, \ldots, p$ are defined in the assumption A4 we have $\sigma > 0$. The convergence of ADMM algorithm for nonconvex function is established in [60] with $\beta = 1$, while we consider $\beta \in (0, 2)$. Our convergence analysis is different and is established under different assumptions. Moreover, we prove the convergence rates which will be discussed in the next section.

## 5 Convergence Rates

In this section, we prove the convergence rates based on the Kurdyka–Łojasiewicz (KŁ) property defined in Sect. 2. We let $u^* := (\mathbf{x}^*, y^*, z^*)$ be the unique limit point of $\{(\mathbf{x}^k, y^k, z^k)\}_{k \geq 0}$, the sequence generated by the PADMM method (5), $v^* := (\mathbf{x}^*, y^*, z^*, y^*, z^*)$, $e_k := \bar{\mathcal{L}}_k - \bar{\mathcal{L}}_*$ and $\lim_{k \to \infty} \bar{\mathcal{L}}_k = \bar{\mathcal{L}}_* := \bar{\mathcal{L}}(v^*)$. We recall readers to see Sect. 2 for definitions on KŁ properties. We begin with the following Lemma.

**Lemma 13** *Suppose that A1–A4 hold, $\bar{\mathcal{L}}$ has the KŁ property at $v^*$ with an exponent $\theta$, and $\sigma$ defined in Lemma 8 is positive. For all $k \geq k_0 \geq 1$ it holds*

$$\bar{\alpha} e_k^{2\theta} \leq e_{k-1} - e_k, \quad \bar{\alpha} = \sigma/(p+2)c^2\tilde{\rho}^2 \tag{42}$$

*where $\tilde{\rho}$ is defined in Lemma 10 and $k_0$ is a large positive integer.*

*Proof* Let $k \geq 1$ fixed. Let $s^k \in \partial \bar{\mathcal{L}}_k$, by Lemma 10 we have

$$\frac{1}{\tilde{\rho}^2(p+2)}|||s^k|||^2 \leq \sum_{i=1}^{p} \|\Delta x_i^k\|^2 + \|\Delta y^k\|^2 + \|\Delta z^k\|^2. \tag{43}$$

By Lemma 8(i) since $\sigma > 0$, we have

$$\sum_{i=1}^{p} \|\Delta x_i^k\|^2 + \|\Delta y^k\|^2 + \|\Delta z^k\|^2 \leq \frac{1}{\sigma}(e_{k-1} - e_k). \tag{44}$$

Combining (43) and (44) leads to

$$\frac{1}{(p+2)\tilde{\rho}^2}|||s^k|||^2 \leq \frac{1}{\sigma}(e_{k-1} - e_k). \tag{45}$$

Let $\epsilon > 0$ arbitrary small. Since $v^k$ converges to $v^*$, there exists a $k_1 \geq 1$ such that $\text{dist}(v^k, v^*) < \epsilon$. By the fact that $\bar{\mathcal{L}}$ has the KŁ property at the point $v^*$ and the sequence $\{\bar{\mathcal{L}}_k\}_{k \geq 1}$ is monotonically decreasing and converges to $\bar{\mathcal{L}}_*$, then there exists $\theta \in [0, 1)$, $c > 0$, and $k_2 \geq 1$ such that for all $k \geq k_2$ it holds

$$(\bar{\mathcal{L}}_k - \bar{\mathcal{L}}_*)^\theta \leq c \, \text{dist}(0, \partial \bar{\mathcal{L}}_k),$$

where it can be equivalently written as $e_k^{2\theta} \leq c^2 \, |||s^k|||^2$. Together with (45), this yields (42) where $k_0 = \max\{k_1, k_2\}$ and $\bar{\alpha} = \sigma/(p+2)c^2\tilde{\rho}^2$. □

**Theorem 2** *Suppose that A1–A4 hold, $\bar{\mathcal{L}}$ has the KŁ property at $v^*$ with an exponent $\theta$, and $\sigma$ defined in Lemma 8 is positive. The sequence $\{e_k\}_{k \geq 1}$ converges to zero in the following rates*

(i) *if $\theta = 0$, then $e_k$ converges to zero within a finite number of iterations;*
(ii) *if $\theta \in (0, 1/2]$, then*

$$e_k \leq e_{k_0}/(1 + \bar{\alpha} e_{k_0}^{2\theta-1})^{k-k_0}, \quad k \geq k_0$$

(iii) *if $\theta \in (1/2, 1)$ then*

$$e_k \leq \left(\mu(k - k_0) + \mu + e_{k_0-1}^{1-2\theta}\right)^{\frac{1}{1-2\theta}}, \quad k \geq k_0$$

*where $\mu > 0$.*

**Proof** (i) Let $\theta = 0$. If $e_k > 0$ for $k \geq k_0$ by Lemma 13 we would have $\bar{\alpha} \leq e_{k-1} - e_k$. As $k$ approaches infinity we have $0 < \bar{\alpha} \leq 0$ which leads to a contradiction. Hence, $e_k = 0$ for all $k \geq k_0$. This implies that there must be a $\tilde{k} \leq k_0$ such that $e_k = 0$ for all $k \geq \tilde{k}$.

(ii) Let $\theta \in (0, \frac{1}{2}]$. The sequence $\{e_k\}_{k\geq 1}$ is monotonically decreasing and since $2\theta - 1 < 0$ then by (42) for all $k \geq k_0$ it holds $\bar{\alpha} e_{k_0}^{2\theta-1} e_k \leq e_k^{2\theta} \leq e_{k-1} - e_k$. We rearrange this to obtain

$$e_k \leq \frac{e_{k-1}}{1 + \bar{\alpha} e_{k_0}^{2\theta-1}} \leq \frac{e_{k-2}}{(1 + \bar{\alpha} e_{k_0}^{2\theta-1})^2} \leq \cdots \leq \frac{e_{k_0}}{(1 + \bar{\alpha} e_{k_0}^{2\theta-1})^{k-k_0}}.$$

This concludes the proof of part (ii).

(iii) Now let $\theta \in (1/2, 1)$. By rearranging (42), we obtain

$$\bar{\alpha} \leq (e_{k-1} - e_k) e_k^{-2\theta}, \quad \forall k \geq k_0 \tag{46}$$

We let $h : \mathbb{R}_+ \to \mathbb{R}$ defined by $h(s) = s^{-2\theta}$ for $s \in \mathbb{R}_+$. Clearly, $h$ is monotonically decreasing as $h'(s) = -2\theta s^{-(1+2\theta)} < 0$ which follows that $h(e_{k-1}) \leq h(e_k)$ for all $k \geq 1$ as $e_k$ is monotonically decreasing. We consider two cases. First, let $r_0 \in (1, +\infty)$ such that $h(e_k) \leq r_0 h(e_{k-1})$, for all $k \geq k_0$. Hence, by (46) we obtain

$$\bar{\alpha} \leq r_0(e_{k-1} - e_k) h(e_{k-1}) \leq r_0 h(e_{k-1}) \int_{e_k}^{e_{k-1}} 1 \, ds$$
$$\leq r_0 \int_{e_k}^{e_{k-1}} h(s) \, ds = r_0 \int_{e_k}^{e_{k-1}} s^{-2\theta} \, ds = \frac{r_0}{1-2\theta} \left( e_{k-1}^{1-2\theta} - e_k^{1-2\theta} \right).$$

Rearrange to get $0 < \frac{\bar{\alpha}(2\theta-1)}{r_0} \leq e_k^{1-2\theta} - e_{k-1}^{1-2\theta}$. Setting $\hat{\mu} = \frac{\bar{\alpha}(2\theta-1)}{r_0} > 0$ and $\nu := 1 - 2\theta < 0$ one then can obtain

$$0 < \hat{\mu} < e_k^\nu - e_{k-1}^\nu, \quad \forall k \geq k_0. \tag{47}$$

Next, we let $h(e_k) \geq r_0 h(e_{k-1})$. This immediately follows that $r_0^{-1} e_{k-1}^{2\theta} \geq e_k^{2\theta}$. We raise both sides to the power $1/2\theta$ and set $q := r_0^{-\frac{1}{2\theta}} \in (0, 1)$ leads to $q e_{k-1} \geq e_k$. Since $\nu = 1 - 2\theta < 0$, then $q^\nu e_{k-1}^\nu \leq e_k^\nu$, and hence,

$$(q^\nu - 1) e_{k-1}^\nu \leq e_k^\nu - e_{k-1}^\nu.$$

By the fact that $q^\nu - 1 > 0$ and $e_p \to 0^+$ as $p \to \infty$, there exists $\bar{\mu}$ such that $(q^\nu - 1) e_{k-1}^\nu > \bar{\mu}$ for all $k \geq k_0$. Therefore, we obtain

$$0 < \bar{\mu} \leq e_k^\nu - e_{k-1}^\nu. \tag{48}$$

Choose $\mu = \min\{\hat{\mu}, \bar{\mu}\} > 0$, one can combine (47) and (48) to obtain

$$0 < \mu \leq e_k^\nu - e_{k-1}^\nu, \quad \forall k \geq k_0.$$

Summing this inequality from $K$ to some $k \geq k_0$ gives $\mu(k - k_0 + 1) + e_{k_0-1}^\nu \leq e_k^\nu$. Thus,

$$e_k \leq (\mu(k - k_0 + 1) + e_{k_0-1}^\nu)^{\frac{1}{\nu}} = \left(\mu(k - k_0 + 1) + e_{k_0-1}^{1-2\theta}\right)^{\frac{1}{1-2\theta}}.$$

This concludes the proof. $\square$

**Theorem 3** *Suppose that A1–A4 hold, and $\bar{\mathcal{L}}$ satisfy the KŁ property at $v^*$ with $\psi \in \Psi_\eta$. If $\sigma > 0$ in Lemma 8, then there exists a $k_0 \geq 2$ such that for all $k \geq k_0$ it holds*

$$|||u^k - u^*||| \leq c \ \max\{\psi(e_{k-1}), \sqrt{e_{k-1}}\}. \tag{49}$$

*Moreover, if $\psi : [0, \eta) \to [0, +\infty)$ is defined by $\psi(s) = s^{1-\theta}$, where $\theta \in [0, 1)$ then the following asymptotic sequential rates are obtained.*

(i) *If $\theta = 0$, then $u^k$ converges to $u^*$ in a finite number of iterations.*
(ii) *If $\theta \in (0, 1/2]$, then*

$$|||u^k - u^*||| \leq e_{k_0} / \left(\sqrt{1 + \bar{\alpha} e_{k_0}^{2\theta-1}}\right)^{k-k_0}, \quad k \geq k_0.$$

*where $\bar{\alpha} = \sigma/(p+2)(1-\theta)^2 \tilde{\rho}^2$.*
(iii) *If $\theta \in (1/2, 1]$, then*

$$|||u^k - u^*||| \leq \left(\mu(k - k_0) + \mu + e_{k_0-1}^{1-2\theta}\right)^{-\frac{1-\theta}{2\theta-1}}, \quad k \geq k_0$$

**Proof** Let $k \geq 1$ be fixed and $\sigma > 0$. By Lemma 8, $\{\bar{\mathcal{L}}_k\}_{k \geq 1}$ is monotonically decreasing and so is $\{e_k\}_{k \geq 1}$, and

$$\frac{1}{p+2}\left(\sum_{i=1}^p \|\Delta x_i^k\| + \|\Delta y^k\| + \|\Delta z^k\|\right)^2 \leq \frac{1}{\sigma}(e_{k-1} - e_k). \tag{50}$$

Let $\epsilon > 0$, $\eta > 0$, and $\psi \in \Psi_\eta$ be given. Since $\bar{\mathcal{L}}_k \searrow \bar{\mathcal{L}}_*$ as $k \to \infty$ and $\lim_{k \to \infty} v^k = v^*$, there exists a $k_0 \in \mathbb{Z}_+$ such that for all $k \geq k_0$

$$\text{dist}(v^k, v^*) < \epsilon \quad \text{and} \quad \bar{\mathcal{L}}_* < \bar{\mathcal{L}}_k < \bar{\mathcal{L}}_* + \eta.$$

By the fact that $\bar{\mathcal{L}}$ satisfies the KŁ property at $v^*$, we then have

$$\psi'(e_k) \cdot \text{dist}(0, \partial \bar{\mathcal{L}}_k) \geq 1. \tag{51}$$

Combining (50) and (51) gives

$$\frac{1}{p+2}\left(\sum_{i=1}^p \|\Delta x_i^{k+1}\| + \|\Delta y^{k+1}\| + \|\Delta z^{k+1}\|\right)^2$$

$$\leq \frac{1}{\sigma}(e_k - e_{k+1})\psi'(e_k) \cdot \mathrm{dist}(0, \partial \bar{\mathcal{L}}_k). \tag{52}$$

Since $\psi \in \Psi_\eta$ is a concave function, then

$$\frac{1}{p+2}\Big(\sum_{i=1}^{p}\|\Delta x_i^{k+1}\| + \|\Delta y^{k+1}\| + \|\Delta z^{k+1}\|\Big)^2$$
$$\leq \frac{1}{\sigma}\big(\psi(e_k) - \psi(e_{k+1})\big) \cdot \mathrm{dist}(0, \partial \bar{\mathcal{L}}_k). \tag{53}$$

For any $\gamma > 0$, by the arithmetic mean-geometric mean inequality this then follows that

$$\sum_{i=1}^{p}\|\Delta x_i^{k+1}\| + \|\Delta y^{k+1}\| + \|\Delta z^{k+1}\|$$
$$\leq \frac{\gamma(p+2)}{2\sigma}\big(\psi(e_k) - \psi(e_{k+1})\big) + \frac{1}{2\gamma}\mathrm{dist}(0, \partial \bar{\mathcal{L}}_k).$$

Thus, by Lemma 10 we then obtain

$$\sum_{i=1}^{p}\|\Delta x_i^{k+1}\| + \|\Delta y^{k+1}\| + \|\Delta z^{k+1}\|$$
$$\leq \frac{\gamma(p+2)}{2\sigma}\big(\psi(e_k) - \psi(e_{k+1})\big)$$
$$+ \frac{\tilde{\rho}}{2\gamma}\Big(\sum_{i=1}^{p}\|\Delta x_i^{k}\| + \|\Delta y^{k}\| + \|\Delta z^{k}\|\Big). \tag{54}$$

Let $\gamma > 0$ large enough such that $\delta := 1 - \frac{\tilde{\rho}}{2\gamma} > 0$. We sum up the latter inequality over $k \geq k_0$ to get

$$\sum_{k=k_0}^{\infty}\Big\{\sum_{i=1}^{p}\|\Delta x_i^{k+1}\| + \|\Delta y^{k+1}\| + \|\Delta z^{k+1}\|\Big\}$$
$$\leq \hat{\delta}\psi(e_{k_0}) + \tilde{\delta}\Big(\sum_{i=1}^{p}\|\Delta x_i^{k_0}\| + \|\Delta y^{k_0}\| + \|\Delta z^{k_0}\|\Big),$$

where $\hat{\delta} := \gamma(p+2)/2\sigma\delta$ and $\tilde{\delta} = \tilde{\rho}/2\gamma\delta$. By the triangle inequality for any $k \geq k_0$, then we have

$$|||u^k - u^*||| \leq \sum_{p \geq k}|||\Delta u^{p+1}|||$$

$$\leq \sum_{p \geq k} \Big\{ \sum_{i=1}^{p} \|\Delta x_i^{p+1}\| + \|\Delta y^{p+1}\| + \|\Delta z^{p+1}\| \Big\}$$

$$\leq \hat{\delta} \psi(e_k) + \tilde{\delta} \Big( \sum_{i=1}^{p} \|\Delta x^k\| + \|\Delta y^k\| + \|\Delta z^k\| \Big).$$

Exploiting (50), the latter inequality leads to

$$|||u^k - u^*||| \leq \hat{\delta} \psi(e_k) + \frac{\tilde{\delta}\sqrt{p+2}}{\sqrt{\sigma}} \sqrt{e_{k-1}} \leq c \max\Big\{\psi(e_k), \sqrt{e_{k-1}}\Big\},$$

where $c := \max\Big\{\hat{\delta}, \frac{\tilde{\delta}\sqrt{p+2}}{\sqrt{\sigma}}\Big\}$. Since $\psi'(s) > 0$ and $\{e_k\}_{k \geq 1}$ is a deceasing sequence, this leads to

$$|||u^k - u^*||| \leq c \max\Big\{\psi(e_{k-1}), \sqrt{e_{k-1}}\Big\}, \quad \forall k \geq k_0.$$

We now let $\psi(s) = s^{1-\theta}$ where $\theta \in [0, 1)$; thus,

$$|||u^k - u^*||| \leq c \max\Big\{e_{k-1}^{1-\theta}, \sqrt{e_{k-1}}\Big\}. \tag{55}$$

Note that $\psi'(s) = (1-\theta)s^{-\theta}$ and since $\bar{\mathcal{L}}$ satisfies the KŁ property at $v^*$, by (51) it holds

$$e_k^{\theta} \leq (1-\theta)|||s^k|||, \quad \forall k \geq k_0$$

where $s^k \in \bar{\mathcal{L}}_k$ and $|||s^k||| = \text{dist}(0, \partial \bar{\mathcal{L}}_k)$. Thus, $e_k^{2\theta} \leq (1-\theta)^2 |||s^k|||^2$ together with Lemma 8 (i) we have $\bar{\alpha} e_k^{2\theta} \leq e_{k-1} - e_k$, where $\bar{\alpha} = \sigma/(p+2)(1-\theta)^2 \tilde{\rho}^2$.

(i) Let $\theta = 0$. Then Theorem 2 (i) gives that $\{e_k\}_{k \geq 1}$ goes to zero in a finite numbers of iterations. Thus, by (55), $u^k$ must converge to $u^*$ in a finite numbers of iterations.

(ii) Let $\theta \in (0, 1/2]$. Note that $\max\{e_{k-1}^{1-\theta}, \sqrt{e_{k-1}}\} = \sqrt{e_{k-1}}$. By Theorem 2(ii),

$$|||u^k - u^*||| \leq \sqrt{e_{k_0}}/\Big(\sqrt{1 + \bar{\alpha} e_{k_0}^{2\theta-1}}\Big)^{k-k_0}, \quad \forall k \geq k_0.$$

(iii) Let $\theta \in (1/2, 1]$, then $\max\{e_k^{1-\theta}, \sqrt{e_{k-1}}\} = e_k^{1-\theta}$. By Theorem 2(iii), we have

$$|||u^k - u^*||| \leq \Big(\mu(k - k_0) + \mu + e_{k_0-1}^{1-2\theta}\Big)^{-\frac{1-\theta}{2\theta-1}}.$$

This completes the proof. □

## 6 Applications

In this section, we consider two nonconvex nonsmooth optimization problems in the form (2) and discuss its associated PADMM algorithm to solve.

**Statistical Learning.** We consider the least square problems with SCAD [25] or MCP [69] regularization functions. The objective function of these problems takes the form

$$\min_x \sum_{i=1}^p r(x_i) + \frac{\mu}{2} \|\mathcal{A}x - y_0\|^2 \qquad (56)$$

where $A \in \mathbb{R}^{m \times n}$ is a measurement matrix, $y_0 \in \mathbb{R}^m$, and $r : \mathbb{R} \to \mathbb{R}$ is a continuous function. For SCAD, the continuous function $r$ takes the form

$$r(t) = r_1(t; \lambda, \theta) = \begin{cases} \lambda|t| & |t| \leq \lambda \\ \frac{-t^2 + 2\theta\lambda|t| - \lambda^2}{2(\theta-1)} & \lambda < |t| \leq \theta\lambda \\ \frac{(\theta+1)\lambda^2}{2} & |t| > \theta\lambda \end{cases}$$

where $\lambda > 0$ and $\theta > 2$, while for MCP, the function $r$ isgiven by

$$r(t) = r_2(t; \lambda, \theta) = \begin{cases} \lambda|t| - \frac{t^2}{2\theta} & |t| \leq \theta\lambda \\ \frac{\theta\lambda^2}{2} & |t| > \theta\lambda \end{cases}$$

where $\lambda > 0$ and $\theta > 0$. By introducing a new variable $y = \mathcal{A}x$, the problem (56) is written as

$$\min_x \sum_{i=1}^p r(x_i) + \frac{\mu}{2} \|y - y_0\|^2, \quad s.t. \quad \mathcal{A}x - y = 0 \qquad (57)$$

Comparing this model with (1), $f(x) = \sum_{i=1}^p r(x_i)$ and $h(y) = \frac{\mu}{2}\|y - y_0\|^2$, both are KŁ functions with exponent $\frac{1}{2}$ [48], $A = \mathcal{A}$, $B = -I_m$, where $I_m \in \mathbb{R}^{m \times m}$ is an identity matrix and $b = 0$. The PADMM algorithm to solve (57) problem with SCAD ($j = 1$) or MCP ($j = 2$) regularization is given by

$$\begin{cases} x^{k+1} = \arg\min_x \sum_{i=1}^p r_j(x_i; \lambda, \theta) + \frac{\alpha}{2}\|\mathcal{A}x - y^k + \alpha^{-1}z^k\|^2 + \frac{1}{2}\|x - x^k\|_Q^2; \\ y^{k+1} = \arg\min_y \frac{\mu}{2}\|y - y_0\|^2 + \frac{\alpha}{2}\|y - \mathcal{A}x^{k+1} - \alpha^{-1}z^k\|^2 + \frac{1}{2}\|y - y^k\|_P^2; \\ z^{k+1} = z^k + \alpha\beta(\mathcal{A}x^{k+1} - y^{k+1}). \end{cases} \qquad (58)$$

We observe that all assumptions A1–A4 are satisfied. The $x$ subproblem involves the second-order term $\frac{\alpha}{2}x^T\mathcal{A}^T\mathcal{A}x$. If $\mathcal{A}^T\mathcal{A}$ is nearly diagonal (or nearly orthogonal) one can replace $\mathcal{A}^T\mathcal{A}$ by a certain symmetric matrix $D$ by setting $Q = \alpha(D - \mathcal{A}^T\mathcal{A})$.

This leads us to an easier computation of $x$. The $y$ subproblem is quadratic and can be solved easily by simply choosing $P = 0$.

**Sparse and low-rank matrix decomposition.** Given $\mathcal{A} \in \mathbb{R}^{m \times n}$ and integer bounds $r > 0$, $s > 0$, the main objective is to find matrices $X_1$, $X_2$ and $Y$ in $\mathbb{R}^{m \times n}$ such that $\mathcal{A} = X_1 + X_2 + Y$, the rank of $X_1$ is less than $r$, the number of nonzero elements in $X_2$ is less than $s$, and change in the columns of $Y$ is relatively small. This problem is formulated as follows

$$\min_{X_1, X_2, Y} \alpha_1 \delta_{\text{rank}(\cdot) \leq r}(X_1) + \alpha_2 \delta_{\|\cdot\|_0 \leq s}(X_2) + \alpha_3 \sum_{i=1}^{n} \|Y_{i+1} - Y_i\|^2,$$
$$\text{s.t.} \quad X_1 + X_2 + Y = \mathcal{A}$$

where $\alpha_i > 0$ for $i = 1, 2, 3$, $\|x\|_0 = \#\{i \mid x_i \neq 0\}$, $Y_i$ is the $i$th column of $Y$. Comparing this model with (1), $f_1(X_1) = \alpha_1 \delta_{\text{rank}(\cdot) \leq r}(X_1)$, $f_2(X_2) = \alpha_2 \delta_{\|\cdot\|_0 \leq s}(X_2)$, and $h(Y) = \alpha_3 \sum_{i=1}^{n} \|Y_{i+1} - Y_i\|^2$ all are KŁ functions with exponent $\frac{1}{2}$. Let $q_1 > 0$ and $q_2 > 0$, then the PADMM algorithm with $Q_1 = q_1 I$, $Q_2 = q_2 I$ and $P = 0$ is given by

$$X_1^{k+1} = \arg\min_{X_1} \alpha_1 \delta_{\text{rank}(\cdot) \leq r}(X_1) + \frac{\alpha}{2} \|X_1 + X_2^k + Y^k + \alpha^{-1} Z^k - \mathcal{A}\|_F^2$$
$$\quad + \frac{q_1}{2} \|X_1 - X_1^k\|^2;$$
$$X_2^{k+1} = \arg\min_{X_2} \alpha_2 \delta_{\|\cdot\|_0 \leq s}(X_2) + \frac{\alpha}{2} \|X_2 + X_1^{k+1} + Y^k + \alpha^{-1} Z^k - \mathcal{A}\|_F^2$$
$$\quad + \frac{q_2}{2} \|X_2 - X_2^k\|^2;$$
$$Y^{k+1} = \arg\min_Y \alpha_3 \sum_{i=1}^{n} \|Y_{i+1} - Y_i\|^2 + \frac{\alpha}{2} \|Y + X_1^{k+1} + X_2^{k+1}$$
$$\quad + \alpha^{-1} Z^k - \mathcal{A}\|_F^2;$$
$$Z^{k+1} = Z^k + \alpha\beta(X_1^{k+1} + X_2^{k+1} + Y^{k+1}).$$

The solution of $X_1$ subproblem is obtained by

$$X_1^{k+1} = \text{proj}_{\{\text{rank}(\cdot) \leq r\}}\left((1 - \alpha\lambda) X_1^k - \alpha\lambda(X_2^k + Y^k + \alpha^{-1} Z^k - \mathcal{A})\right),$$

where $0 < \lambda < 1/(\alpha + q_1)$ and the projection onto $\{rank(\cdot) \leq r\}$ is done by the Singular Value Decomposition. The solution of $X_2$ subproblem is obtained by

$$X_2^{k+1} = \text{proj}_{\{\|\cdot\|_0 \leq s\}}\left((1 - \alpha\gamma) X_2^k - \alpha\gamma(X_1^{k+1} + Y^k + \alpha^{-1} Z^k - \mathcal{A})\right),$$

where $0 < \gamma < 1/(\alpha + q_2)$. To project onto $\{\|\cdot\|_0 \leq s\}$ we can set all the elements except the $s$ largest ones (in absolute value) equal to zero. The $Y$ subproblem is simple quadratic so we omit discussion on how to solve it.

## 7 Concluding Remarks

We have considered a multi-block proximal ADMM algorithm for solving linearly constrained separable nonconvex nonsmooth optimization problems. We established the global convergence analysis by showing that the sequence generated by PADMM is Cauchy. The convergence rate results for the function and the sequence are provided in the KŁ framework. We observe that the convergence rates are directly related to the KŁ exponent of the modified augmented Lagrangian. In the future, we plan to develop some calculus rules for deducing the KŁ exponent for the augmented Lagrangian in a nonconvex setting. The latter will allow us to determine the convergence rates for the PADMM.